\documentclass{elsart}
\usepackage{amsmath, amssymb, latexsym, amsfonts,graphicx}
\DeclareGraphicsExtensions{.eps, .ps}
\reversemarginpar
\newcommand{\up}[1]{\raisebox{.5ex}{{\scriptsize #1}}}

\newcommand{\tr}{\mbox{}^t\!}

\newcommand{\R}{\mathbb R} 
\newcommand{\C}{\mathbb C} 
\newcommand{\Z}{\mathbb Z} 

\newcommand{\Rnn}{{\mathbb R}_{\geqslant 0}}
\newcommand{\RP}{\mathbb R P} 
\newcommand{\CP}{\mathbb C P}
\newcommand{\Ztwo}{(\mathbb Z_2)}

\newcommand{\bd}{\partial}
\newcommand{\sig}{\sigma}
\newcommand{\gT}{\mathbf{T}}
\newcommand{\cT}{\mathcal{T}}

\newcommand{\one}{\textbf{(1)}}
\newcommand{\two}{\textbf{(2)}}

\newcommand{\inc}{\subset}
\newcommand{\sign}{\mathop{\mathrm{sign}}}

\begin{document}

\begin{frontmatter}
\title{Ambient Surfaces and T-Fillings \\
              of T-Curves}
\author{Bertrand Haas}
\address{Mathematical Science Research Institute,\\ 
1000 Centennial Dr., Berkeley, CA 94720, USA}

\begin{abstract}
\emph{T-curves} are piecewise linear curves which have been used with success
since the beginning of the 1990's to construct new real algebraic curves
with prescribed topology mainly on the real projective plane
(see~\cite{Shustin/curv-prescr-sing,Itenberg/Ragsdale1,Haas/multiluc}).  In
fact T-curves can be used on any real projective toric surface.  We
generalize here the construction of the latter by departing from non-convex
polygons and we get \emph{ambient surfaces} that we characterize
topologically and on which T-curves are also well defined.  Associated to
each T-curve we construct a surface, the \emph{T-filling}, which is
topologically the analog of the quotient of the complexification of a real
algebraic curve by the complex conjugation.  We use then the T-filling
\textbf{(1)} to prove a theorem for T-curves on arbitrary ambient surfaces
which is the analog of a theorem of Harnack for algebraic curves, and
\textbf{(2)} to define the orientability and orientation of a T-curve, the
same way it is defined in the real algebraic setting.
\end{abstract}
\end{frontmatter}

\noindent\textbf{\large Table of content:} \par\noindent
Introduction \\
Section \ref{part1}:  The Ambient Surfaces \\
\null\hspace{4mm} Subsection \ref{def-latt-geom} (p. \pageref{def-latt-geom}):  
    Some Definitions in Lattice Geometry \\
\null\hspace{4mm} Subsection \ref{loc-topo-amb} (p. \pageref{loc-topo-amb}):  
    The Ambient Surface and its Local Topology \\
\null\hspace{4mm} Subsection \ref{glob-topo-amb} (p. \pageref{glob-topo-amb}):  
    The Global Topology of an Ambient Surface \\
Section \ref{part2}:  The T-Curves  \\
\null\hspace{4mm} Subsection \ref{latT-amb} (p. \pageref{latT-amb}):  
    T-Curves on an Ambient Surface \\
\null\hspace{4mm} Subsection \ref{isom-Tcurv} (p. \pageref{isom-Tcurv}):  
    Isomorphic T-curves \\
\null\hspace{4mm} Subsection \ref{Tcurv-RP2} (p. \pageref{Tcurv-RP2}):  
    T-curves on $ \RP^2 $ \\
Section \ref{F(K)}:  The T-Filling of a T-curve  \\
\null\hspace{4mm} Subsection \ref{const-fill} (p. \pageref{const-fill}):  
    The Construction \\
\null\hspace{4mm} Subsection \ref{F(K)-alg} (p. \pageref{F(K)-alg}):  
    The Relation with Algebraic Geometry \\
\null\hspace{4mm} Subsection \ref{Tfill-applic} (p. \pageref{sec:Har-T}):  
    Applications

\vspace{5mm}
\noindent\textbf{\Large Introduction}
\vspace{5mm}

In the $ 1980 $'s Oleg Viro invented a geometric method (called
\emph{patchworking}) to construct smooth algebraic hypersurfaces of real
projective toric varieties~\cite{Viro/these2} by gluing, like in a patchwork
quilt, different toric varieties with given smooth embedded hypersurfaces.
In particular, his method, when used with lines as patches, gives rise to a
family of piecewise linear curves which have been commonly called
\emph{T-curves}.  T-curves and related constructions have been used with
success since the beginning of the $ 1990 $'s to answer various questions
around the first part of Hilbert 16\up{th} Problem~\cite{Hilbert/pbs} (more
precisely about the topology of real algebraic curves and surfaces, see for
instance~\cite{Shustin/curv-prescr-sing,Itenberg/Ragsdale1,Haas/multiluc,%
Itenberg/T-surfaces}, and most recently about the number of limit cycles in
polynomial vector fields (see~\cite{Shustin-Itenberg/lim-cycl-polyn-vf}).

T-curves originally lie on (polygonal models of) real projective toric
surfaces.  They are constructed from the data of \textbf{(1)} an integral
convex polygon, \textbf{(2)} a rectilinear triangulation of the polygon where
the vertices of the triangles are the integer points of the polygon, and
which satisfies some convexity property (see def.~\ref{conv-trig}), and
\textbf{(3)} a distribution of signs (``$+$'' and ``$-$'') on the integer
points of the polygon.  The convex polygon must be viewed at the same time as
a polygon defining the ambient toric surface and as the Newton polygon of
some polynomial $ \sum a_{i,j} x^i y^j $ (the Newton polygon is defined as
the convex hull of the integer points $ (i,j) $ such that $ a_{i,j} \not = 0
$).  The signs on the integer points $ (i,j) $ must be viewed as the signs of
the $ a_{i,j} $'s.  The polynomial defines an algebraic curve on the toric
surface, which is \emph{congruent} to the T-curve, i.e.\ there is a
homeomorphism of the toric surface which transforms the T-curve into the
algebraic curve.  Notice that not all algebraic curves are congruent to
T-curves with same Newton polygon.

In the first part of his 16\up{th} problem Hilbert asked what are all the
possible non-singular real algebraic curves of given degree on the real
projective plane $ \RP^2 $ up to congruence.  This problem generalizes to
curves with various constraints on various surfaces (for instance curves with
a given Newton polygon on a given real projective toric surface).  We study
here, in this spirit, T-curves on some combinatorial surfaces, the
ambient surfaces (see below), which generalize some classical models of
non-singular real projective toric surfaces.

In the construction of a T-curve the restrictions to convex polygons and
convex triangulations are necessary for the parallel construction of the
polynomial defining the algebraic curve congruent to the T-curve.  But these
two restrictions are not necessary for the construction of the T-curve
itself.  Jesus De Loera~\cite{Deloera/convex-patchw} has studied some
properties of T-curves on $ \RP^2 $ without the assumption of convexity on
the triangulation.  The purpose of this approach is not any more the
construction of new algebraic curves, but the proof of properties of
T-curves, in particular the ones inspired directly from algebraic curves on $
\RP^2 $ (for instance De Loera proved an analog of Hilbert's lemma). 

The same approach would be interesting for T-curves on arbitrary real
projective toric surfaces, but since we dropped the link with algebraic
curves by dropping the convexity requirement on the triangulation, it is also
natural to drop the link with real projective toric surfaces by dropping the
requirement that the Newton polygon be convex.  This is what we do in the
first part of this article and we get this way the \emph{ambient surfaces} on
which T-curves are still well defined.  

In this more general setting we construct in part II a surface from a
T-curve, the \emph{T-filling}, which is the analog of the quotient of the
complexification of an algebraic curve by the complex conjugation.  The
latter surface is well known and has been classically used to prove many
properties of real algebraic curves.  Here we use the T-filling to prove in a
similar way the analog of a theorem of Harnack~\cite{Harnack1} for algebraic
curves on $ \RP^2 $: The maximum number of connected components of a T-curve
on an arbitrary ambient surface is equal to the number of integer points in
the interior of its Newton polygon (see subsect.~\ref{sec:Har-T}).
For each Newton polygon a remarkable T-curve, called the \emph{Harnack
T-curve} (see~\ref{Har-const}), has the maximum number of connected
components (such a curve is called a \emph{maximal curve}).  Harnack T-curves
are well known on $ \RP^2 $ where they are congruent to curves constructed by
Harnack~\cite{Harnack1}.  We also use the T-filling to define
\emph{orientability} and \emph{orientations} of T-curves (see
subsect.~\ref{sec:orient-T}) in the same way it has been defined by
Rokhlin~\cite{Rokhlin/cplex-orient-curv} for algebraic curves on $ \RP^2 $.

Viro's patchworking method for curves, because of its real algebraic setting,
deals with convex polygons.  Allowing polygons to be non-convex leads to a
patchworking method for T-curves with more freedom.  In particular we show
in~\cite{Haas/PhD} that every maximal T-curve is the patchworking of Harnack
T-curves (most of which have non-convex Newton polygons).  This is a
first step in characterizing which ones among all maximal real algebraic
curves on a given projective toric surfaces can be constructed as T-curves
(indeed this \emph{Harnack patchworking} holds in particular for maximal
T-curves with convex Newton polygons and with convex triangulations).

The T-filling is a new object.  The Harnack theorem for T-curves on $ \RP^2 $
and the notion of orientability and orientation for T-curves on $ \RP^2 $
were known by I.~Itenberg and O.~Viro but, as far as I know, never published
by them (Jes\'us De Loera discusses orientability and orientation of T-curves
on $ \RP^2 $ after Itenberg and Viro in~\cite{Deloera/convex-patchw}).
\pagebreak

\section{The Ambient Surfaces}
\label{part1}

\subsection{Some Definitions in Lattice Geometry}
\label{def-latt-geom}

\begin{defn}
Here all segments are \emph{integral rectilinear segments}, that is their two
endpoints are integer points.  A \emph{(closed) polygonal line} is a finite
union of segments $ s_1 \cup \dots \cup s_r \inc \R^2 $, such that for $ i =
2, \dots, r-1 $, (for $ i $ modulo $ r $), the segment $ s_i $ shares one of
its endpoints with $ s_{i-1} $ and its other endpoint with $ s_{i+1} $.  We
always assume that the union of any two consecutive segments $ s_i \cup
s_{i+1} $ is not a segment.  Here we call a \emph{polygon} a subset of $ \R^2
$ bounded by a closed polygonal line and homeomorphic to a disk.  An
\emph{edge of the polygon} is a segment of its boundary polygonal line.  A
\emph{vertex of a polygon} is an endpoint of an edge.
\end{defn}

Notice that the vertices of a polygon are all integer points. 

\begin{defn}
A \emph{triangulation of a polygon} is a cell decomposition
of the polygon into simplexes.  A simplex of dimension $ 2 $ is a
\emph{triangle of the triangulation}, a simplex of dimension $ 1 $ is an
\emph{edge of the triangulation}, and a simplex of dimension $ 0 $ is a
\emph{vertex of the triangulation}.
\end{defn}

Vertices and edges of a triangulation shouldn't be mixed up with vertices and
edges of a polygon.  We are going now to define a \emph{parity} which assigns
to various objects an element of $ \Ztwo^2 = (\Z / 2\Z) \times (\Z / 2\Z) $.
We denote the four values of the parities by $ (0,0) $, $ (0,1) $, $ (1,1) $
and $ (1,0) $.  A parity is said to be \emph{even} if its value is $ (0,0) $.
A parity is \emph{odd} if it is not even.

\begin{defn}
The \emph{parity of an integer point} $ (p_1, p_2) $ is $ (p_1 \bmod 2, p_2
\bmod 2) $.  The \emph{parity of a segment} is the sum of the two values that
the parity takes on all the integer points of the segment.
\end{defn} 

Notice that the parity of the integer points which lie on a given segment
takes exactly two values.  Also, notice that the parity of a segment is never
even.

\begin{defn}
The \emph{parity of a vertex of a polygon} is the sum of the parities of the
two segments underlying the edges of the polygon meeting in the vertex.  The
\emph{parity of an edge of a polygon} is the sum of the parities of the two
vertices of the polygon which end the edge.
\end{defn}

Notice that the parity of a vertex (of an edge) of a polygon can possibly
take all four values.  Also the parity of a vertex (of an edge) of a polygon
may be different than the parity of its underlying integer point (of its
underlying segment).

\begin{lem}
\label{nev-rpar-2nev}
If a polygon contains a vertex of odd parity, then it contains at least two
vertices of odd parity.
\end{lem}

\begin{pf}
Let us follow the boundary of the polygon, departing from the vertex
of odd parity and coming back to it.  Since the vertex has odd parity, the
"departing-edge" has different parity than the "coming-back-edge".  So on the
way we must pass through a vertex where two edges of different parity meet,
i.e. through a vertex of odd parity.
\qed
\end{pf}

\begin{defn}
The vertices of odd parity of a polygon divide its boundary into
connected components, the closure of which we call the \emph{broken edges of
the polygon}.  
\end{defn}

Notice that the parity takes the same value on all the segments of a given
broken edge.

\begin{defn}
\label{rpar-brok}
The \emph{parity of a broken edge of a polygon} is the sum of the parities
of all the edges of the polygon contained in the broken edge.
\end{defn}

Notice that if a polygon has vertices of odd parity, then the parity of a
broken edge of the polygon is equal to the sum of the parities of the two
vertices of the polygon which end the broken edge.

\subsection{The Ambient Surface and its Local Topology}
\label{loc-topo-amb}

\subsubsection{The Ambient Surface Associated to a Polygon}
\label{intro-amb}

Let $ \Pi $ be a polygon in $ \Rnn^2 := \{(x,y) \in \R^2 | x\geq 0, y\geq 0\}
$, and for every $ a,b \in \{0,1\} $, let $ \sig_{a,b} $ be the symmetry $
(x,y) \mapsto ((-1)^a x, (-1)^b y) $.  We glue the disjoint union of the four
copies $ (\sig_{a,b} \cdot \Pi) $ by their boundary in the following way:  we
identify each point $ (x,y) $ of the boundary with the point $ \sig_{(a,b)}
\cdot (x,y) $ where $ (a,b) $ is such that $ \langle (a,b), (c,d) \rangle :=
ac + bd = 0 $, when $ (x,y) $ lies on an edge of parity $ (c,d) $.  We obtain
in this way a surface without boundary.

\begin{defn}
\label{df:amb-surf}
We call \emph{ambient surface} (associated to a polygon $ \Pi $), the
surface constructed above, and we denote it by $ S(\Pi) $.  
\end{defn}

The four maps $ \sig_{a,b} $ glue to a map $ \mu: S(\Pi) \rightarrow \Pi $
which is a fourfold ramified covering.  The ramifications take place along
the broken edges in the following way:  A point in the interior of a broken
edge lifts to two points, and an endpoint of a broken edge lifts to one
point.

\begin{center}
\begin{picture}(4.6,2)
\put(0,1.6){$ \bigsqcup_{a,b} (\sig_{a,b} \cdot \Pi) $}
\put(2.5,1.7){\vector(1,0){1.1}}
\put(2.9,1.85){$ \scriptstyle q $}
\put(3.8,1.6){$ S(\Pi) $}
\put(1.5,1.4){\vector(2,-1){2.3}}
\put(2.7,.9){$ \scriptstyle \pi $}
\put(4.15,1.4){\vector(0,-1){1}}
\put(4.3,.9){$ \scriptstyle \mu $}
\put(4,0){$ \Pi $}
\end{picture}
\end{center}

\begin{figure}[h]
\includegraphics{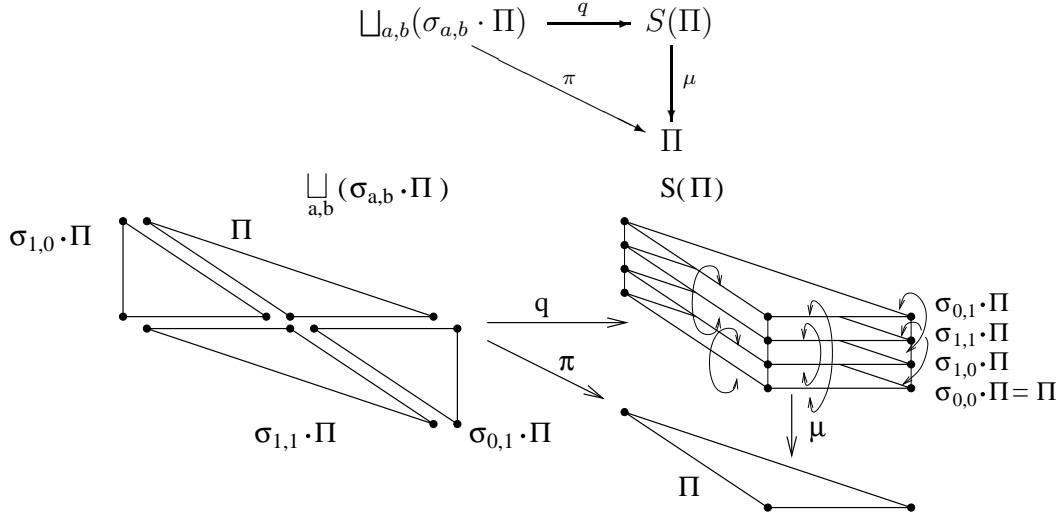}
\caption{\small The fourfold ramified covering structure of $ S(\Pi) $.}
\label{figure1}
\end{figure}

\begin{figure}[h]
\includegraphics{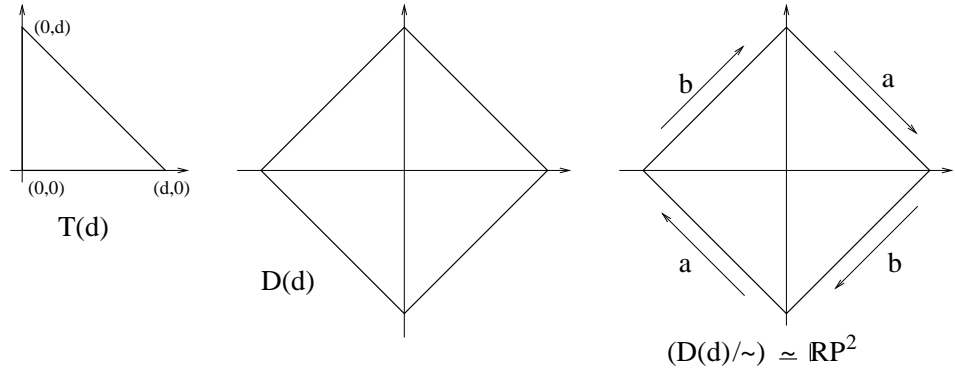}  
\caption{\small Recovering the real projective plane $ \RP^2 $ from the
triangle $ T_d $.}
\label{figure2}  
\end{figure}

\begin{rem}
\label{amb-tor}
If the polygon $ \Pi $ is convex and the broken edges are just the edges of $
\Pi $, then the ambient surface $ S(\Pi) $ is a topological model of the
real projective toric surface $ X(\Pi) $  whose image by the moment map is $
\Pi $.  If $ \Pi $ is convex but some broken edge is the union of at least
two edges of $ \Pi $, then $ X(\Pi) $ contains some conic points and $ S(\Pi)
$ is a desingularization of $ X(\Pi) $.  The following example will be used
later.  Let $ d $ be a positive integer and $ \gT_d $ be the triangle in $
\R^2 $ with vertices $ (0,0) $, $ (d,0) $ and $ (0,d) $ then it is well known
that $ X(\gT_d) = \RP^2 $.  Let's check that $ S(\gT_d) $ is homeomorphic to 
$ \RP^2 $:  Let $ D_d $ be the \emph{diamond} in $ \R^2 $ equal to the union
of the four copies $ (\sig_{a,b} \cdot T_d) $~, $ a,b \in \{0,1\} $.  Then $
S(T_d) $ is obtained from $ D_d $ by identifying every two opposite points on
its boundary (see fig.~\ref{figure2}).  As another example let $
\widetilde{\gT}_d $ be the triangle in $ \R^2 $ with vertices $ (0,0) $, $
(2d,0) $ and $ (0,d) $, then $ X(\widetilde{\gT}_d) $ is a pinched torus.
Indeed the chart at $ (0,d) $ gives an affine cone, and the charts at $ (0,0)
$ and $ (2d,0) $ give affine planes.  It is easy to check that $ S(\Pi) $ is a
sphere.  The topology of an arbitrary ambient surface is described in
proposition~\ref{char-amb}.
\end{rem}

\begin{defn}
The image of a copy $ (\sig_{a,b} \cdot \Pi) $ by the quotient map $ q :
\bigsqcup_{a,b} (\sig_{a,b} \cdot \Pi) \rightarrow S(\Pi) $ is a
\emph{quadrant} of $ S(\Pi) $.  Since $ q $ is one-to-one from a given copy $
(\sig_{a,b} \cdot \Pi) $ to its image, we will identify (abusively) a
quadrant with its pre-image $ (\sig_{a,b} \cdot \Pi) $.  
\end{defn}

\subsubsection{The Local Structure around the Lift of a Vertex}
\label{loc-vert}

Let $ u $ be a vertex of $ \Pi $, and $ \ell $ and $ \ell' $ be the two edges
of $ \Pi $ meeting at $ u $.  

\begin{figure}[h]
\includegraphics{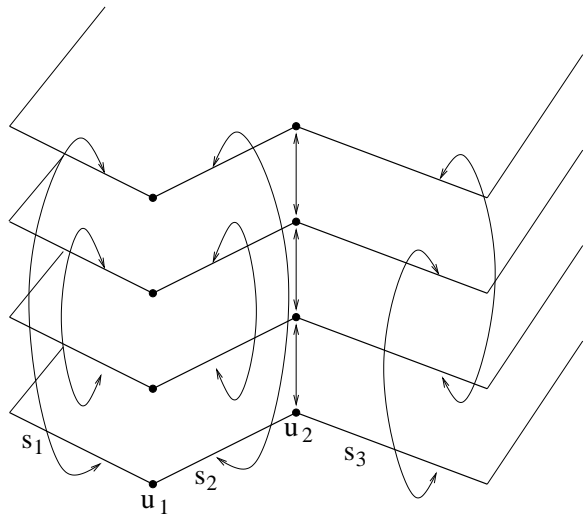}
\caption{\small How the segments in the lift of a broken edge are identified
two by two according to their parity.}
\label{figure3}
\end{figure}

Assume first that $ u $ is of even parity, so $ \ell $ and $ \ell' $ are of same
parity.  The construction of $ S(\Pi) $ (see def.~\ref{df:amb-surf}) implies
then that two copies of $ \Pi $ will be glued to one another by identifying
the two corresponding copies of $ \ell \cup \ell' $, and the two other copies
of $ \Pi $ will be glued to one another by identifying the two other copies
of $ \ell \cup \ell' $ (see fig.~\ref{figure3}).

Assume now that $ u $ is of odd parity.  The construction of $ S(\Pi) $ (see
def.~\ref{df:amb-surf}) implies that the four copies of $ \Pi $ are glued to
one another in the following way:  The union of $ \ell $ with a reflection of
$ \ell $ through a coordinate axis is identified to the union of the two
other copies of $ \ell $ and the union of $ \ell' $ with a reflection of $
\ell' $ through the other coordinate axis is identified to the union of the
two other copies of $ \ell' $ (see fig.~\ref{figure3}).

\begin{figure}[h]
\includegraphics{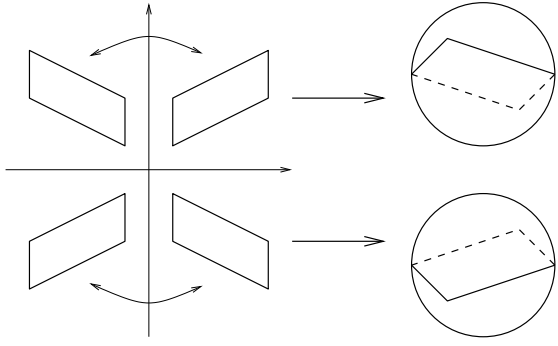}
\caption{\small Gluing the four copies of $ \Pi $ to one another gives a
union of two spheres.}
\label{figure4}
\end{figure}

If all the segments of $ \bd \Pi $ are of the same parity (that is if $ \Pi $
has only one broken edge), then two copies of $ \Pi $ will be glued to one
another by all their edges, and the two other copies of $ \Pi $ will be also
glued to one another by all their edges. We get this way two connected
components, each homeomorphic to a sphere (see an example on
fig.~\ref{figure4}).  

Assume now that the parity of the segments of $ \bd \Pi $ takes at least two
values.  Notice that since the two endpoints $ u $ and $ u' $ of a broken
edge $ \ell $ are vertices of odd parity, they lift each to only one point in $
S(\Pi) $.  All the other points of $ \ell $ lift to two points.  Therefore the
lift $ \mu^{-1}(\ell) $ of $ \ell $ in $ S(\Pi) $ is a circle which is a
ramified twofold covering of $ \ell $, where the ramifications take place in
the lifts $ \mu^{-1}(u) $ and $ \mu^{-1}(u') $.

\subsubsection{Canonical Charts for the Ambient Surface}
\label{charts}

Let $ \ell_1, \dots, \ell_r $ be the broken edges of $ \Pi $, indexed such
that $ \ell_i $ meets $ \ell_{i+1} $~, $ i \bmod r $, and let $ u_i $ be the
vertex of $ \Pi $ where they meet.  Let
$$ 
U_i = S(\Pi) \setminus \big(\bigcup_{\stackrel{\scriptscriptstyle j \not = i}%
{\scriptscriptstyle j \not = i+1}}
\mu^{-1}(\ell_j)\big) 
$$ 
So $ U_i $ is an open neighborhood of $ \mu^{-1}(u_i) $ in $ S(\Pi) $ and is
homeomorphic to $ \R^2 $.  We will assume that the broken edges are oriented
locally near their endpoints, away from their endpoints.  So $
\mu^{-1}(\ell_i) $ and $ \mu^{-1}(\ell_{i+1}) $ are oriented in $ U_i $.

\begin{defn}
The open neighborhood $ U_i $ equipped with a homeomorphism $ U_i
\rightarrow \R^2 $ mapping $ \mu^{-1}(\ell_i) $ and $ \mu^{-1}(\ell_{i+1}) $, with
their orientations, onto the coordinate axis $ Ox $ and $ Oy $ is a
\emph{(oriented) chart of $ S(\Pi) $}.
The system of charts $ U_1, \dots, U_r $ will be called a \emph{canonical
system of charts}
for $ S(\Pi) $.  
\end{defn}

Notice that the chart $ U_i $ tells how to glue the four quadrants of $
S(\Pi) $ onto $ \mu^{-1}(\ell_i) $ and $ \mu^{-1}(\ell_{i+1}) $.

\begin{defn}
We call \emph{(open) quadrants of the chart} $ U_i
$ and denote $ U_i^{0,0} $, $ U_i^{0,1} $, $ U_i^{1,1} $, and $ U_i^{1,0} $
the connected components of $ U_i \setminus (\mu^{-1}(\ell_i) \cup
\mu^{-1}(\ell_{i+1})) $ such that $ U_i^{a,b} $ is mapped to the set $ \{(x,y)
\in \R^2, (-1)^a x > 0, (-1)^b y > 0\} $ by the homeomorphism $ U_i
\rightarrow \R^2 $.  
\end{defn}

Notice that the closure in $ S(\Pi) $ of an open quadrant $ U_i^{a,b} $ is
a quadrant $ (\sig_{c,d} \cdot \Pi) $ of $ S(\Pi) $.  So the chart is
determined by the correspondence $ (c,d) \mapsto (a,b) $.

\begin{figure}[h]
\includegraphics{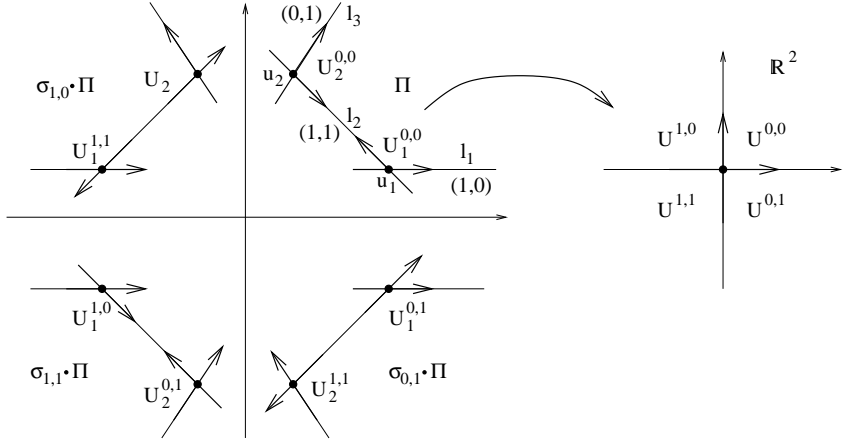}
\caption{\small A local orientation of all the broken edges determines the
label of the quadrants by the $ U_i^{a,b} $.  The parity of the broken
edges $ \ell_i $ is in parentheses.}
\label{figure5}
\end{figure}

\begin{defn}
We call \emph{the parity matrix of the chart} $ U_i $,
the matrix 
$ \renewcommand{\arraystretch}{1}
M_i = \Big(
\begin{array}{cc}
\alpha_i & \alpha_{i+1} \\
\beta_i  & \beta_{i+1}
\end{array}
\Big) 
$
where $ (\alpha_j,\beta_j) $ is the parity of the segments of the broken edge
$ \ell_j $.  \end{defn}

\begin{lem}
\label{quad-chart-amb}
If the closed quadrant $ \bar{U}_i^{a,b} $ of the chart $ U_i $ is equal to
the quadrant $ (\sig_{c,d} \cdot \Pi) $ of $ S(\Pi) $, then $ (a,b) = ((c,d)
\cdot M_i) $.
\end{lem}

\begin{pf}
With our choice of local orientations for the broken edges, $
\bar{U}_i^{0,0} $ is always equal to $ \Pi = (\sig_{0,0} \cdot \Pi) $, so the
permutation $ (c,d) \mapsto (a,b) $ can be considered as an element of $
GL(2, \Z_2) $.  Since $ \bar{U}_i^{0,1} $ and $ \bar{U}_i^{1,0} $ correspond
to the quadrants $ (\sig_{c,d} \cdot \Pi) $ glued to quadrant $ \Pi $ along $
\mu^{-1}(\ell)_i $ and $ \mu^{-1}(\ell)_{i+1} $ respectively, we get from the
construction of $ S(\Pi) $ (see def.~\ref{df:amb-surf}) that 
$$ \renewcommand{\arraystretch}{1.2}
(c,d) = \bigg\{ \begin{array}{cc}
		   (\beta_i, \alpha_i)     & \,\text{if}\, (a,b) = (0,1) \\
		   (\beta_{i+1}, \alpha_{i+1}) & \,\text{if}\, (a,b) = (1,0)
		   \end{array} 
$$
which is equivalent to $ (c,d) = (a,b) \cdot M_i^{-1} $.
\qed
\end{pf}

Notice that the transformations $ M_i \mapsto M_j $ correspond to the
gluing of chart $ U_i $ with chart $ U_j $.

\begin{defn}
We will call $ M_{i,j} : M_i \mapsto M_j $ the \emph{gluing
transformation} of the parity matrices $ M_i $
and $ M_j $.
\end{defn}

\subsubsection{Gluing Charts of the Ambient Surface}
\label{glu-charts}

Assume now that $ \Pi $ has more than two broken edges.  Let $ U_1, \dots,
U_r $ be a canonical system of charts of $ S(\Pi) $.  

Then the gluing of two consecutive charts $ U_i $ and $ U_{i+1} $ is easily
stated: 

\begin{lem}
\label{gl-2charts}
The gluing transformation $ M_{i,i-1} $ is given by the matrix 
$ \renewcommand{\arraystretch}{.8}
 \Big(   \begin{array}{cc}
	    0 & 1 \\
	    1 & \eta_i 
	    \end{array} \Big)
$ and the right-product:  $ M_i = M_{i-1} \cdot M_{i,i-1} $, where $ \eta_i =
0 $ if the broken edge $ \ell_i $ has even parity, and $ \eta_i = 1 $ if $ \ell_i $
has odd parity.
\end{lem}

\begin{pf}
We assume without loss of generality that $ i = 2 $, so we will
compute with:
$$ 
\renewcommand{\arraystretch}{1}
M_1 = \Big( \begin{array}{cc}
		\alpha_1 & \alpha_2 \\
		\beta_1  & \beta_2 
		\end{array} \Big) \qquad\text{and}\qquad 
M_2 = \Big( \begin{array}{cc}
		\alpha_2 & \alpha_3 \\
		\beta_2  & \beta_3
		\end{array} \Big)
$$
Notice that the parities $ (\alpha_1, \beta_1) $ and $ (\alpha_3, \beta_3) $
of the segments of the broken edges $ \ell_1 $ and $ \ell_3 $ are equal if and only
if $ \ell_2 $ has even parity, so the lemma holds when $ \eta_2 = 0 $.  Assume
that $ (\alpha_3, \beta_3) \not = (\alpha_1, \beta_1) $, so $ \eta_2 = 1 $.
Since the parities belong to $ \Ztwo^2 $ and are never even, and since $
(\alpha_1, \beta_1) \not = (\alpha_2, \beta_2) $, then $ (\alpha_3, \beta_3)
= (\alpha_1, \beta_1) + (\alpha_2, \beta_2) $.  
\qed
\end{pf}

\begin{figure}[h]
\includegraphics{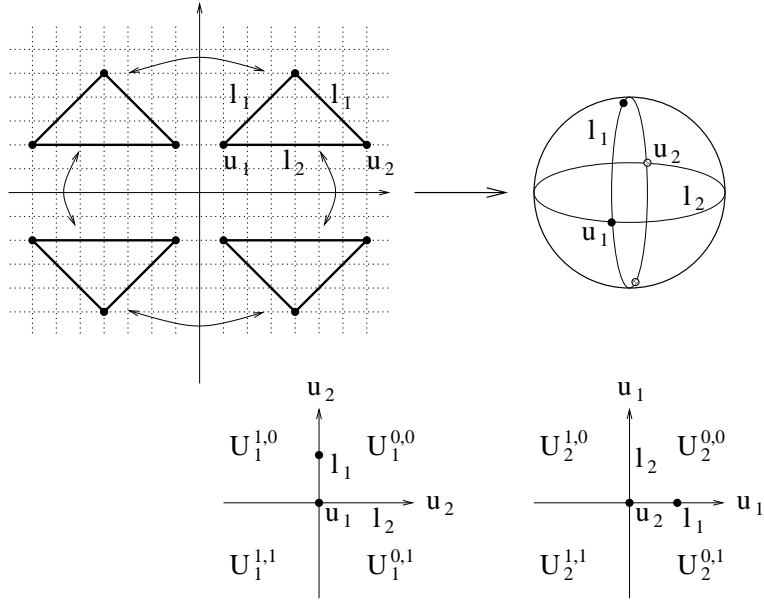}
\caption{\small Gluing the four copies of $ \Pi $ to one another gives a
sphere.}
\label{figure6}
\end{figure}

Notice that if a polygon $ \Pi $ has only two broken edges, then
the canonical system of charts for $ S(\Pi) $ has only two charts.  Since
the parity of the segments of $ \bd \Pi $ takes only two values the gluing of
the two charts gives a surface $ S(\Pi) $ which is a sphere (see
fig.~\ref{figure6}).

Notice that the gluing transformation of any parity matrices is equal to
the product of the transformations of consecutive parity matrices:
$$
M_{i,j} = \prod_{k=0}^{k=j-i-1} M_{i+k, i+k+1} 
$$

\begin{defn}
The canonical system of charts $ U_i $ equipped with the gluing
transformations $ M_{i,j} $ is the \emph{canonical atlas
structure} of $ S(\Pi) $.  
\end{defn}

\begin{cor}
\label{rpar-tubli}
If $ \ell_i $ has even parity, then a tubular neighborhood of
$ \mu^{-1}(\ell_i) $ is an annulus. \\
If $ \ell_i $ has odd parity, then a tubular neighborhood of
$ \mu^{-1}(\ell_i) $ is a Moebius band.
\end{cor}

\begin{pf}
With our choice of local orientations for the broken edges of $ \Pi $
we have always 
$$ U_i^{0,0} = U_{i-1}^{0,0} \quad\text{and}\quad U_i^{0,1} = U_{i-1}^{1,0}
$$
and either 
$$ U_i^{1,0} = U_{i-1}^{0,1} \quad\text{and}\quad U_i^{1,1} = U_{i-1}^{1,1}
$$ 
In which case the tubular neighborhood of $ \ell_i $ is an annulus (see
fig.~\ref{figure7}), or
$$ U_i^{1,0} = U_{i-1}^{1,1} \quad\text{and}\quad U_i^{1,1} = U_{i-1}^{0,1}
$$
In which case the tubular neighborhood of $ \ell_i $ is a Moebius band (see
fig.~\ref{figure7}).

\begin{figure}[h]
\includegraphics{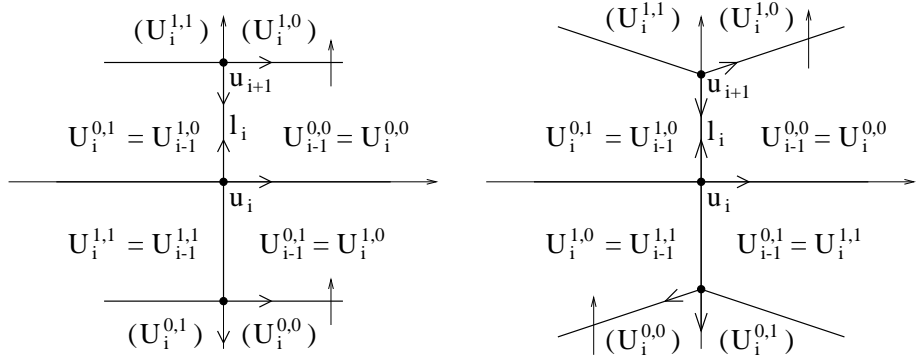} 
\caption{\small The tubular neighborhood of the lift $ \mu^{-1}(\ell) $ of a
broken edge $ \ell $ is either an annulus if $ \ell $ is even, either a
Moebius band if $ \ell $ is odd.}
\label{figure7} 
\end{figure}

The closed quadrants $ \bar{U}^{0,0} $ are always represented by $ \Pi $.
Let $ (\sig_{c,d} \cdot \Pi) $ be the quadrant representing $
\bar{U}_{i-1}^{0,1} $, so $ (0,1) = (c,d) \cdot M_{i-1} $.  Now $ (\sig_{c,d}
\cdot \Pi) $ represents also a closed quadrant $ \bar{U}_i^{a,b} $, so $
(a,b) = (c,d) \cdot M_i $.  Since $ M_i = M_{i-1} \cdot M_{i,i-1} $ we get
that $ (a,b) = (0,1) \cdot M_{i,i-1} $.  From lemma~\ref{gl-2charts} we get
directly that $ (a,b) = (1,0) $ if $ \ell_i $ is even (hence we get an
annulus), and $ (a,b) = (1,1) $ if $ \ell_i $ is odd (hence we get a Moebius
band).
\qed
\end{pf}

Notice that, given a canonical system of charts, the gluing transformations
of consecutive charts are determined by the $ 0 $--$ 1 $ sequence of the
numbers $ \eta_i $ of lemma~\ref{gl-2charts}.  So two polygons with the same
sequence of $ \eta_i $'s up to circular permutation give rise to homeomorphic
ambient surfaces.

\subsection{The Global Topology of an Ambient Surface}
\label{glob-topo-amb}

\subsubsection{A Basis for the 1-Homology of the Ambient Surface}
\label{1-hom-basis}

Let $ S(\Pi) $ be an ambient surface, and let $ U_1, \dots, U_r $ be a
canonical system of charts on $ S(\Pi) $.  We have seen that when $ r = 1 $,
then $ S(\Pi) $ is a homeomorphic to the disjoint union of two
spheres~(\ref{figure4}), and when $ r = 2 $, then $ S(\Pi) $ is homeomorphic
to a single sphere (fig.~\ref{figure6}).  Assume now that $ r \geq 3 $.  Recall
that with the notations of the charts, the $ \ell_i $'s are the broken edges of
$ \Pi $ with endpoints $ u_{i-1} $ and $ u_i $.

\begin{lem}
\label{basis-homol}
The lifts $ \mu^{-1}(\ell_3), \dots, \mu^{-1}(\ell_r) $ of the broken edges $ \ell_3,
\dots, \ell_r $ form a basis for the 1-homology space $ H_1(S(\Pi)) $ (with
coefficients in $ \Z $ if $ S(\Pi) $ is orientable, and in $ \Z_2 $ if $
S(\Pi) $ is not orientable).
\end{lem}

\begin{pf}
To prove this lemma it suffices to check two points:
\begin{itemize}

\item  The set $ \{\mu^{-1}(\ell_3), \dots, \mu^{-1}(\ell_r)\} $ is free in $
	H_1 $ since $ S(\Pi) \setminus (\cup_3^r \mu^{-1}(\ell_i)) $ is
	connected.  

\item  The set $ \{\mu^{-1}(\ell_3), \dots, \mu^{-1}(\ell_r)\} $ is complete
	since $ S(\Pi) \setminus (\cup_3^r \mu^{-1}(\ell_i)) $ is simply
	connected.

\end{itemize}
\qed
\end{pf}

\subsubsection{The Topology of Ambient Surfaces}
\label{topo-amb}

Assume here that $ \Pi $ has at least two broken edges (recall that if it has
only one, $ S(\Pi) $ is homeomorphic to the disjoint union of two spheres).

\begin{lem}
\label{2brokodd} 
\begin{itemize}

\item  If $ \Pi $ has a broken edge of odd parity, then it has two
	broken edges of odd parity.

\item  If $ \Pi $ has two consecutive broken edges of odd parity, then it
	has three broken edges of odd parity.

\end{itemize}
\end{lem}

\begin{pf}
Let $ U_1, \dots, U_r $ be a canonical system of charts on the
surface $ S(\Pi) $.  From lemma~\ref{gl-2charts}, we get that, to each broken
edge of odd parity, correspond a gluing matrix 
\renewcommand{\arraystretch}{.9}
$ A_1 = \Big( \begin{array}{cc}
		 0 & 1 \\
		 1 & 1
		 \end{array} \Big) $ 
and to each broken edge of even parity
correspond a gluing matrix 
$ A_0 = \Big( \begin{array}{cc}
		 0 & 1 \\
		 1 & 0
		 \end{array} \Big) $.
The product $ \prod $ of all gluing matrices of consecutive charts must be
the identity.  The product with $ A_0 $ just permutes the columns.  So if
$ \prod $ has a matrix $ A_1 $ as factor then it must have at least two
matrices $ A_1 $ as factors.  Since 
$ A_1^2 = \Big( \begin{array}{cc}
		 1 & 1 \\
		 1 & 0
		 \end{array} \Big) $,
if $ \prod $ has two matrices $ A_1 $ as consecutive factors, then it must
have at least three matrices $ A_1 $ as factors.
\qed
\end{pf}

\begin{prop}
\label{char-amb}
Let $ r > 1 $ be the number of broken edges of $ \Pi $.
\begin{itemize}

\item  The surface $ S(\Pi ) $ is orientable if and only if the parity of
	the segments of $ \bd \Pi $ takes only two values.   

\item  If $ S(\Pi) $ is orientable, then $ r $ is even, and $ S(\Pi) $ is
	the connected sum of $ r/2 - 1 $ tori.  

\item  If $ S(\Pi) $ is not orientable, then it is the connected sum of $ r
	-2 $ projective planes.

\end{itemize}
\end{prop}

\begin{pf}
From corollary~\ref{rpar-tubli} we know that the tubular
neighborhoods of the lifts of the broken edges of $ \Pi $ are all annuli if
and only if all the broken edges are of even parity.  From
lemma~\ref{basis-homol}, we know that the lifts of all the broken edges but
two consecutive ones is a basis of the 1-dimensional homology space of $
S(\Pi) $.  From lemma~\ref{2brokodd} we know that if all but two consecutive
broken edges are known to be of even parity, then all the broken edges are of
even parity.  

Therefore $ S(\Pi) $ is orientable if and only if all its broken edges are
of even parity.  From definition~\ref{rpar-brok} it is clear that
the broken edges are all of even parity if and only if the parity
of the segments of $ \bd \Pi $ takes only two values. 

From lemma~\ref{basis-homol} we know that the dimension of the
1-dimensional homology space of $ S(\Pi) $ is $ r-2 $.  Since $ r > 1 $, the
parity of the segments of $ \bd \Pi $ takes more than one value, and then $
S(\Pi) $ is connected.  So if $ S(\Pi) $ is orientable, it is the connected
sum of $ (r-2)/2 $ tori, and if it is non-orientable, it is the connected sum
of $ r-2 $ projective planes.
\qed
\end{pf}

\begin{rem}
\label{generaliz-amb}
Notice that in the construction of an ambient surface, the requirement that
the polygon be homeomorphic to a disk can be dropped.  We could also
construct ambient surfaces from generalized polygons, the boundary of which
would be a collection of integral closed polygonal lines.  However we tried
to avoid here a too wide generalization.  This is motivated by the necessity
of this requirement for the ambient surfaces we use in~\cite{Haas/PhD} to
construct maximal T-curves.  The construction of ambient surfaces can also be 
generalized to higher dimensional \emph{ambient varieties} by using piecewise
integral linear subsets of $ \R^n $.  We may require certain constraints on
these subsets.  For instance we could ask that only $ n $ faces of dimensions
$ n-1 $ meet at each vertex (so we keep a nice local topological structure),
or that each $ n $-face be homeomorphic to a $ n $-ball, etc.
\end{rem}

\section{The T-Curves}
\label{part2}

\subsection{T-Curves on an Ambient Surface}
\label{latT-amb}

\begin{defn}
A \emph{primitive triangulation} of a polygon $ \Pi $ is a
triangulation of the polygon such that the set of vertices of the
triangulation is equal to the set of integer points $ \Pi \cap \Z^2 $.
\end{defn}

\begin{figure}[ht]
\includegraphics{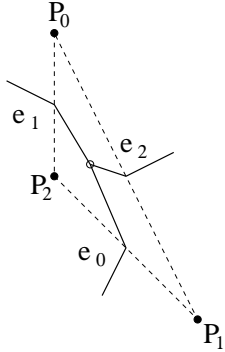}  
\caption{$ \sign(e_i) = \delta(P_j) \delta(P_k) $, so $ \prod \sign(e_i) =
\prod \delta(P_i)^2 = 1 $.  Therefore either $ 0 $ or $ 2 $ edges have
negative sign.}
\label{figure8}  
\end{figure}

Choose a polygon $ \Pi \inc \Rnn^2 $, and let $ S(\Pi) $ be our ambient
surface. We now define a T-curve on $ S(\Pi) $ in five steps: 

\begin{itemize}
\item  Choose a primitive triangulation $ \cT $ of $ \Pi $.  The reflections
    through the two coordinate axis generate a triangulation of the disjoint
    union $ \bigsqcup_{a,b} (\sig_{a,b} \cdot \Pi) $~, $ a,b \in
    \{0,1\} $, which induces a triangulation of $ S(\Pi) $.  With the
    notations of section~\ref{intro-amb} the triangulation of $
    \bigsqcup_{a,b} (\sig_{a,b} \cdot \Pi) $ is $ \pi^{-1}(\cT) $, and the
    triangulation of $ S(\Pi) $ is $ \mu^{-1}(\cT) $.

\item  Choose a sign function $ \Pi \cap \Z^2 \rightarrow \{+1,-1\} \; ,\;
    (x,y) \mapsto \delta(x,y) $.  We extend this distribution on $
    \bigsqcup_{a,b} (\sig_{a,b} \cdot \Pi) \cap \Z^2 $ in the following way:
    \begin{eqnarray} 
    \label{ext-sd2}
    \hfill \delta((-1)^a x, (-1)^b y) 
    & = & (-1)^{ax + by} \delta(x,y) \cr
    \text{Or equivalently:}\qquad \delta(\sig_{a,b} \cdot (x,y)) 
    & = & (-1)^{\langle (a,b), (c,d) \rangle} \delta(x,y) 
    \end{eqnarray}
    Where $ a,b \in \{0,1\} $~, $ (c,d) $ is the parity of the point $ (x,y)
    $, and $ \langle (a,b), (c,d) \rangle = ac + bd $.

\item  We assign to each edge of the triangulation $ \pi^{-1}(\cT) $  the
    sign equal to the product of the signs of its endpoints.  This induces a
    sign function on the edges $ \mu^{-1}(\cT) $.  Indeed for an arbitrary
    edge $ e $ of $ \pi^{-1}(\cT) $ which lies on the boundary of $
    \pi^{-1}(\Pi) $, we check that the sign of the edge $ \sig_{c,d} $ to
    which $ e $ becomes identified in $ \mu^{-1}(\cT) $ is the same than the
    sign of $ e $.  From the construction of $ S(\Pi) $ (see
    def.~\ref{df:amb-surf}) we get that the parity $ (a,b) $ of $ e $ and $
    (c,d) $ verify $ \langle (a,b), (c,d) \rangle = 0 $, and from
    formula~\ref{ext-sd2} above we get that 
    \begin{equation}
    \label{e-sign}
    \sign(\sig_{c,d} \cdot e) = (-1)^{\langle (a,b), (c,d) \rangle} \sign(e)
    \end{equation}
    which is then equal to $ \sign(e) $.

\item The sign function on the edges of $ \mu^{-1}(\cT) $ has the property
    that each triangle has either $ 0 $ or $ 2 $ edges of negative sign.
    Indeed the product of the three signs of the edges of a triangle is equal
    to $ +1 $ since it is the product of the squares of the three signs of
    the vertices of the triangle (see fig.~\ref{figure8}).  

\item Consider a cell decomposition of $ S(\Pi) $ dual to the triangulation $
    \mu^{-1}(\cT) $, and let $ G $ be the graph underlying its 1-skeleton.
    So each edge of $ G $ is dual to an edge of $ \mu^{-1}(\cT) $.  We assign
    to the former edge the sign of the latter edge.  From the previous step
    we deduce that each vertex of $ G $ is the endpoint of exactly two edges
    of negative sign and one edge of positive sign.  Therefore the union of
    the edges of negative sign form a closed curve $ K $ on $ G $ (and hence
    on $ S(\Pi) $).  
    
\end{itemize}

Notice that the data $ (\Pi, \cT, \delta) $ and $ (\Pi, \cT, -\delta) $
define the same curves.  

\begin{defn}
\label{defT}
The \emph{T-curve} defined by the data $ (\Pi, \cT, \pm\delta) $ is the curve
constructed as above.  It will be denoted by $ K(\Pi, \cT, \delta) $, (or $
K(\Pi, \cT, -\delta) $). 
\end{defn}

\subsubsection{T-Curves and Algebraic Curves}
\label{Tcurv-alg-curv}

\begin{defn}
\label{conv-trig}
A triangulation of a convex polygon $ \Pi $ is said to be \emph{convex}, or
\emph{regular}, if there exist a convex function $ \nu : \Pi \rightarrow \R $
which is linear on each triangle, and non linear on any union of two
triangle.  
\end{defn}

\begin{rem}
\label{Tcurv/tor}
If $ \Pi $ is a convex polygon in $ \R^2 $, then we have seen
in remark~\ref{amb-tor}  that the ambient surface $ S(\Pi) $ is a topological
model of (i.e. is homeomorphic to) the real projective toric surface $ X(\Pi)
$ defined by $ \Pi $.  In fact the construction of T-curves was initially
motivated by the construction of new real algebraic curves on real projective
toric surfaces as shows \emph{Viro's theorem}: Let $ K = K(\Pi, \cT, \delta)
$ be a T-curve.  If $ \cT $ is a convex triangulation, then $ K $ is
congruent to a real algebraic curve with Newton polygon $ \Pi $.  For
instance T-curves constructed from a convex triangulation of the triangle $
\gT_d $ (see remark~\ref{amb-tor}) are congruent to real algebraic curves of
degree $ d $.
\end{rem}

\begin{defn}
\label{T-Newt-pol}
In the remark above, observe that the convex polygon $ \Pi $ is
the Newton polygon of an algebraic curve.  By analogy, we will call the
polygon $ \Pi $ in the data defining a T-curve $ K = K(\Pi, \cT, \delta) $
the \emph{Newton polygon} of $ K $.
\end{defn}

\subsubsection{A digression to higher dimensions}
\label{sec:high-dim}

As we can generalize ambient surfaces to higher dimensional ambient
varieties (see remark~\ref{generaliz-amb}), we can also generalize the
construction of T-curves to \emph{T-hypersurfaces}.  The data defining a
T-hypersurface is similar:  A Newton polytope $ \Pi $, a triangulation $
\cT $ of $ \Pi $ and a distribution of signs on the integer points of $ \Pi
$.  The triangulation $ \cT $ is a cell decomposition of $ \Pi $ into
simplexes where the vertices of the simplexes are all the integer points
of the Newton polytope.  Notice that simplexes of dimension more than two may
have different volumes.  We can add more constraints and require that all
simplexes of highest dimension of the triangulation have same volume.  

Notice that Viro's theorem (see remark~\ref{Tcurv/tor}) for T-curves requires
that the convex function (see def.~\ref{conv-trig}) over the Newton polygon
have integer values on integer points.  This is not a serious restriction.
Indeed we can allow rational values on the integer points since there are
finitely many integer points and we can multiply the function by the least
common multiple of the denominators.  A perturbation of the vertices of a
three dimensional polytope does not destroy the incidence structure of the
faces if it is small enough.  Therefore we can also allow the function to
take real values on the integer points since an arbitrary small perturbation
of it will transform the real values into rational values.

The definition of convex triangulation of a convex polytope in dimension more
than two is similar than in dimension two (see def.~\ref{conv-trig}): There
should exist a convex real valued function on the Newton polytope which is
affine on each simplex and non affine on any union of two simplexes.  Notice
that in high dimensions it is not true that the incidence structure of the
faces of a polytope is preserved under a small enough perturbation of the
vertices.  Therefore it is safer to require in the definition of a convex
triangulation in higher dimension that the convex function takes rational
values on the vertices.

The connection of T-hypersurfaces with real algebraic hypersurfaces on
projective toric varieties relies again on the convexity of the triangulation
of the Newton polytope.  In fact Viro's theorem is stated in all dimensions:
If the triangulation of the Newton polytope of the T-hypersurface is convex,
then there exist a real algebraic hypersurface which is congruent to the
T-hypersurface.  

\subsection{Isomorphic T-curves}
\label{isom-Tcurv}

\subsubsection{Translation of the Newton polygon}
\label{translat-Newt}

Let $ \Pi $ be a polygon in $ \Rnn^2 $, and $ (s,t) $ be an
integral vector such that the translated polygon $ \Pi' = \Pi + (s,t) $ lies
also in $ \Rnn^2 $, and let $ K(\Pi, \cT, \delta) $ be a T-curve.  Let $ \cT'
$ be the triangulation of $ \Pi' $ translated from $ \cT $, and let $ \delta'
$ be the sign distribution on $ \Pi' \cap \Z^2 $ defined by $ \delta'(x,y) =
\delta(x-s, y-t) $.  

\begin{prop}
\label{translat-Pi}
The T-curves $ K = K(\Pi, \cT, \delta) $ and $ K' = K(\Pi', \cT', \delta') $
are equal.
\end{prop}

\begin{pf}
It is clear that the construction of the ambient surface $ S(\Pi) $
depends on $ \Pi $ up to translation by an integral vector, so $ S(\Pi') =
S(\Pi) $.  Let $ (c,d) $ be the parity of the point $ (s,t) $, and let $ a,b
\in \{0,1\} $.  If $ (x,y) \in \Pi $, then 
$$ 
\delta'((-1)^a (x+s), (-1)^b (y+t)) = (-1)^{\langle (a,b), (c,d) \rangle}
\delta((-1)^a x, (-1)^b y) 
$$
Recall that the sign of an edge of a triangulation is the product of the
signs of its two endpoints.  Therefore the edges of the triangulation of $
S(\Pi') $ have the same sign than the corresponding edges of the
triangulation of $ S(\Pi) $, so $ K = K' $.
\qed
\end{pf}

\subsubsection{Linear Transformations of the Newton Polygon}
\label{lin-transf-Newt}

If $ K $ is a T-curve with Newton polygon $ \Pi $, and $ c,d \in
\{0,1\} $, let $ K^{c,d} $ be the restriction of $ K $ to the quadrant $
(\sig_{c,d} \cdot \Pi) $.  

Let $ K(\Pi, \cT, \delta) $ be a T-curve and let $ A \in GL(2,\Z)
$.  We will note $ A_2 $ the reduction of $ A $ in $ GL(2,\Z_2) $.  Let $
 \cT' $ be the triangulation of $ (A \cdot \Pi) $
resulting from applying $ A $ to $ \cT $, let $ \delta'(x,y) = \delta(A^{-1}
\cdot (x,y)) $, and let $ A \cdot K $ be the T-curve defined by $ (A \cdot
\Pi, \cT', \delta') $.  

\begin{prop}
\begin{itemize}

\item  $ A \cdot K $ and $ K $ are congruent.  

\item  The homeomorphism $ S(\Pi) \rightarrow S(A \cdot \Pi) $ which
transforms $ K $ to $ A \cdot K $, transforms each restriction $ K^{c,d} $
into a restriction $ (A \cdot K)^{s,t} $. 

\item  $ (c,d) = (s,t) \cdot A_2 $.

\end{itemize}
\end{prop}

\begin{pf}
The action of $ A_2 $ on parities underlies the action of $ A $, so
an edge $ e $ of the triangulation $ \cT $ is of even parity if and only if
the edge $ A \cdot e $ of $ \cT' $ is of even parity.  So it is clear from
prop.~\ref{char-amb} that $ S(A \cdot \Pi) $ and $ S(\Pi) $ are homeomorphic.

Let $ U $ and $ U' $ be canonical charts of $ S(\Pi) $ and of $ S(A \cdot \Pi)
$ around an arbitrary vertex (of odd parity) $ u $ and around $ A \cdot u $
respectively.  For any $ a,b \in \{0,1\} $, let $ (\sig_{c,d} \cdot \Pi) $
be the quadrant of $ S(\Pi) $ representing $ \bar{U}^{a,b} $, and $
(\sig_{s,t} \cdot (A \cdot \Pi)) $ be the quadrant of $ S(A \cdot \Pi) $
representing $ (\bar{U}')^{a,b} $.  Let $ M $ be the parity matrix of $ U $.
Then $ (A_2 \cdot M) $ is the parity matrix of $ U' $.  From
lemma~\ref{quad-chart-amb} we get that 
$$ 
(a,b) = (c,d) \cdot M \quad\text{and}\quad (a,b) = (s,t) \cdot (A_2 \cdot M) 
$$
Therefore $ (c,d) = (s,t) \cdot A_2 $.

To prove the proposition, we must show that the T-curve $ K $ is congruent
to $ A \cdot K $ by a homeomorphism $ U \rightarrow U' $ such that for any
$ a,b \in \{0,1\} $~, $ U^{a,b} \mapsto {U'}^{a,b} $.  It suffices then to
show that the sign of an edge of the triangulation of $ (\sig_{c,d} \cdot
\Pi) $ is equal to the sign of the corresponding edge of the triangulation of
$ (\sig_{s,t} \cdot (A \cdot \Pi)) $:  

Let $ e $ be a segment of $ \bd \Pi $ and $ (a,b) $ its parity.  So $ ((a,b)
\cdot \tr A_2) $ is the parity of $ (A \cdot e) $.  From
the definition of the sign distribution $ \delta' $ of $ (A \cdot K) $, we get
that the sign of $ (A \cdot e) $ is equal to the sign of $ e $.  Then from
the construction in section~\ref{intro-amb} we get that

\begin{eqnarray*} 
\sign(\sig_{s,t} \cdot (A \cdot e)) 
& = & (-1)^{\langle (s,t) , (a,b) \cdot \tr A_2 \rangle} \sign(A \cdot e)  \\
& = & (-1)^{\langle (s,t) \cdot A_2, (a,b) \rangle} \sign(e) \\
& = & (-1)^{\langle (c,d), (a,b) \rangle} \sign(e) \\
& = & \sign(\sig_{c,d} \cdot e)
\end{eqnarray*}
\qed
\end{pf}

\subsection{T-Curves on $ \RP^2 $}
\label{Tcurv-RP2}

\begin{defn}
\label{df:prim-seg}
A \emph{primitive segment} is a segment such that the only integer points it
contains are its two endpoints.  The \emph{integral length} (or simply
\emph{length}) of a polygonal line is the number of primitive segments
contained in it. Two integer points are \emph{neighbors} if they can be
connected by a primitive segment.  
\end{defn}

\begin{defn}
A piece of a curve homeomorphic to the segment $ [0,1] \inc \R $ is called
an \emph{arc} of the curve.
\end{defn}

\begin{defn}
\label{df:ovals}
Let $ K $ be a collection of disjoint embedded circles in a
surface.  A connected component of $ K $ is called an \emph{oval}
if it bounds a subset of the surface homeomorphic to a disk.  The interior of
this subset is called the \emph{inside} of the
oval.  The interior of the complementary of the inside is the
\emph{outside} of the oval.  An oval
with no other ovals inside will be called an \emph{empty
oval}, and an oval inside no other ovals will be called an
\emph{outermost oval}.  A connected component of $ K $
which is not an oval is called a \emph{nontrivial component} of $ K $.  
\end{defn}

Recall from remark~\ref{amb-tor} that $ \gT_d $ is the triangle in $ \R^2 $
with vertices $ (0,0) $, $ (d,0) $ and $ (0,d) $, and that $ S(\gT_d) $ is
a model of the projective plane.

\begin{prop}
Let $ K(\gT_d, \cT, \delta) $ be a T-curve on $ \RP^2 $.  Then
\begin{itemize}

\item  $ K $ has only ovals as connected components if its degree
$ d $ is even.

\item  $ K $ has one and only one nontrivial component if its degree $ d $
is odd.

\end{itemize}
\end{prop}

\begin{pf}
It is clear that on $ \RP^2 $ two nontrivial embedded circles which
intersect transversally, intersect an odd number of time, and an oval which
is intersected by any embedded circle transversally is intersected an even
number of times.  Therefore, since $ K $ is a disjoint union of embedded
circles, it has at most one nontrivial connected component.

Consider now a simple path on the diamond $ D(d) $ through the edges of the
triangulation $ \cT $, going from one point of the boundary of the diamond
to its opposite point.  This path lifts to a loop on $ RP^2 $ which is a
nontrivial embedded circle and if it cuts $ K $, it cuts it transversally.
So the loop intersects $ K $ an odd number of times if and only if $ K $ has
a nontrivial connected component.

Let us follow the path from one endpoint to the other.  If the loop
intersects the T-curve $ r $ times, then the integer points on the path
change signs $ r $ times.  So $ r $ is odd if and only if the two
endpoints of the path have opposite signs.  The formula~\ref{ext-sd2}
page~\pageref{ext-sd2} for the extension of the sign distribution from the
triangle $ \gT(d) $ to the diamond $ D(d) $ shows that two opposite points on
the boundary of $ D(d) $ have opposite signs if an only if the degree $ d $
of the T-curve is odd.
\qed
\end{pf}

Notice that this lemma holds if we replace ``T-curve'' by `` nonsingular
algebraic curve''.  In that case it can be proved in a similar way, which is
a direct consequence of Bezout theorem applied to the curve and a generic
line.

Notice now that if an oval $ O $ of a T-curve lies inside a quadrant of the
ambient surface then all the integer points lying inside $ O $ and outside
any oval which lies inside $ O $ have same sign.  Therefore the following
definition is coherent.

\begin{defn}
\label{df:sign-ov}
Let $ O $ be an oval of a T-curve, which lies inside a quadrant of the
ambient surface.  The \emph{sign of the oval} $ O $ is the sign of the
integer points lying inside $ O $ and outside any oval which lies inside $ O
$.
\end{defn}

\section{The T-Filling of a T-Curve.}
\label{F(K)}

Recall from definition~\ref{defT} that a T-curve $ K(\Pi, \cT, \delta) $ lies
on the graph $ G $ which underlies the 1-skeleton of a dual cell
decomposition of the triangulation $ \cT_S $ of the ambient surface $ S =
S(\Pi) $.  We note here $ \cT_S $ for the triangulation of $ S $ (that is $
\cT_S = \mu^{-1}(\cT) $ where $ \mu $ is the projection $ S \rightarrow \Pi
$).  The dual cell decomposition we use here is the one which can be refined
into the barycentric cell decomposition of $ \cT_S $.  Therefore an edge of $
G $ is the union of two segments linking the barycenters of two adjacent
triangles $ t $ and $ t' $ and passing through the middle of the edge $ t
\cap t' $ of $ \cT_S $.

\begin{defn}
\label{df:incidgr-trig}
The \emph{incidence graph of the triangulation} $ \cT_S $ is the graph
refined from $ G $ by subdividing the edges of $ G $ at the points of $ G
\cap \cT_S $, (so the edges are now segments).  We note this graph $ G(S,
\cT_S) $ or simply\footnote{The notation is a little abusive since $ G $ is
also a graph on $ S $ but acceptable since $ G $ and $ G(S) $ have same
geometric realization (the unions of their edges are equal).} $ G(S) $.  The
image by the projection $ \mu $ of $ G(S) $ will be called the
\emph{incidence graph of the triangulation} $ \cT $ and will be noted $
G(\Pi, \cT) $ or simply $ G(\Pi) $.  Hence $ G(S) $ is equal to the lift $
\mu^{-1}(G(\Pi)) $.
\end{defn}

In accordance with section~\ref{latT-amb}, every edge of $ G(S)
$ inherit of the sign of the supporting edge of $ G $ (that is the sign of
the edge of $ \cT $ it intersects).  Notice that each edge of $ G(\Pi) $
lifts to four edges of $ G(S) $.

\begin{lem}
\label{2-e-neg}
For each edge of $ G(\Pi) $, there is exactly two edges, among its four lifts
in $ G(S) $, which have negative signs.
\end{lem}

\begin{pf}
Let $ e $ be an edge of $ \cT $, and let $ (a,b) $ be its parity.
Recall the formula~\ref{e-sign} page~\pageref{e-sign} which gives the sign of
a copy $ (\sig_{c,d} \cdot e) $ of $ e $:
$$ 
\sign(\sig_{c,d} \cdot e) = (-1)^{\langle (a,b), (c,d) \rangle}
\sign(e) 
$$ 
Notice that $ \langle (a,b), (c,d) \rangle = 0 $ if and only if $ (c,d) =
(0,0) $ or $ (c,d) = (b, a) $.  Since the parity of a segment is never even,
$ (b, a) \not = (0,0) $.  Therefore there is exactly two edges, among the
four which lift from $ e $, which have same sign than $ e $ (so there is also
two edges which have opposite sign than $ e $).  This proves the assertion,
since every edge of $ G(\Pi) $ inherits of the sign of the edge
of $ \cT $ that it intersects.
\qed
\end{pf}

\subsection{The Construction} 
\label{const-fill}

Let $ K = K(\Pi, \cT, \delta) $ be a T-curve and $ S = S(\Pi) $ its ambient
surface.  Each component $ \tilde{\gamma} $ of $ K $ is a cycle (in the language
of graph theory) in the incidence graph $ G(S) $, i.e. $ \tilde{\gamma} $ is
a cyclic sequence of edges of $ G(\Pi) $, each edge sharing one
end with the previous edge and the other end with the next edge.  These
cycles don't intersect.  The projection $ \gamma = \mu(\tilde{\gamma}) $ on $
\Pi $ will be considered as the cycle made up by the projections of the edges
of $ \tilde{\gamma} $.  These cycles may intersect. 

\begin{defn}
a \emph{thick-Y} is a tubular neighborhood of the three edges of $
G(\Pi) $ lying in a triangle of the triangulation $ \cT $.
\end{defn}

More precisely let $ e_1, e_2 $ and $ e_3 $ be the three edges of $ G(\Pi) $
lying in a triangle $ t $ of $ \cT $, and indexed counterclockwise, and let $
A $ and $ B_i $ be the endpoints of $ e_i $.  The thick-Y is obtained by
gluing the three ribbons $ r_i = e_i \times [-1,+1] $ by the following way:
identify $ A \times [-1,0] $ in $ r_i $ with $ A \times [0,+1] $ in $ r_{i+1}
$ by $ (A,x) \mapsto (A,-x) $ (see fig.~\ref{figure10}).

\begin{figure}[ht]
\includegraphics{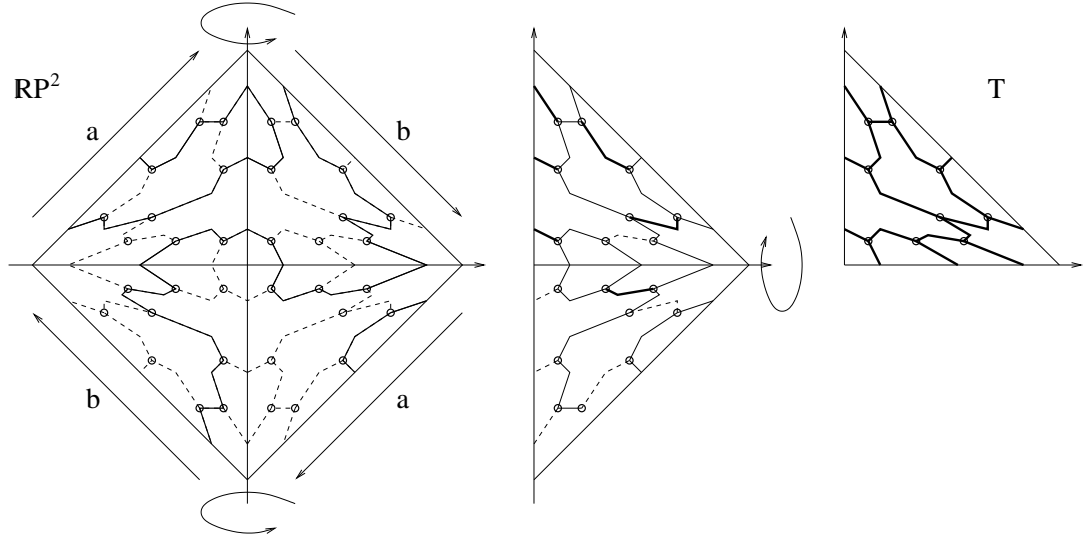}  
\caption{\small The T-curve (plain lines) on the incidence graph $ G(S) $
(dotted lines).  The folded T-curve on $ G(\gT_d) $ uses every edge twice
(thick plain lines).}
\label{figure9}  
\end{figure}

\begin{figure}[h]
\includegraphics{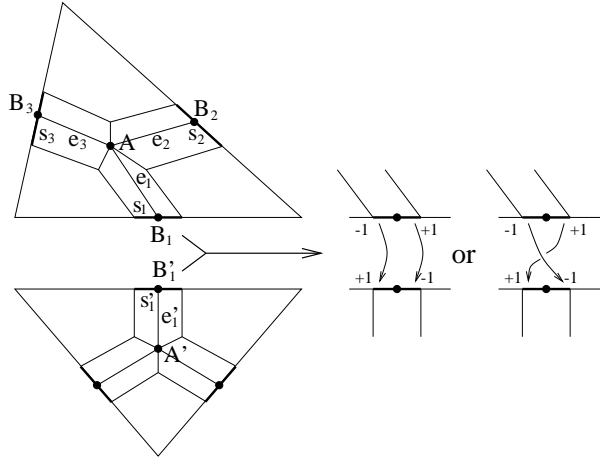}  
\caption{\small A thick-Y is glued from three ribbons, and two thick-Y's are
glued with or without a twist.}
\label{figure10}  
\end{figure}

From lemma~\ref{2-e-neg} we get that every edge of $ G(\Pi) $ is used
twice by the cycles $ \gamma $ (see fig.~\ref{figure9}).  Then every union of
two edges of $ G(\Pi) $ which lie in a triangle is an arc of a cycle $ \gamma
$.  In other words every cycle $ \gamma $ is the union of such arcs.  We
now describe how to glue adjacent thick-Y's together according to which arcs
are consecutive on the cycles.

Let $ t' $ be a triangle of $ \cT $ adjacent to $ t $, and let $ e'_1, e'_2 $
and $ e'_3 $ be the three edges of $ G(\Pi) $ lying in $ t' $, indexed
counterclockwise, with endpoints $ A' $ and $ B'_i $, and such that $ B'_1 =
B_1 $.  Let $ s_1 = B_1 \times [-1,+1] $ be an end-segment of the thick-Y in
$ t $, and let $ s'_1 = B'_1 \times [-1,+1] $ be the corresponding
end-segment of the thick-Y in $ t' $.  From the disjoint union of the two
thick-Y's in $ t \cup t' $, we identify $ s_1 $ with $ s'_1 $: 
\begin{enumerate}

\item  either by $ (B_1,x) \mapsto (B'_1,-x) $ if $ (e_2 \cup e_1) $ and $
	e'_1 \cup e'_3) $ are two consecutive arcs of some cycle~$ \gamma $.
	So $ (e_3 \cup e_1) $ and $ (e'_1 \cup e'_2) $ are two consecutive
	arcs of some cycle $ \gamma' $.

\item  either by $ (B_1,x) \mapsto (B'_1,x) $ if $ (e_2 \cup e_1) $ and $
	(e'_1 \cup e'_2) $ are two consecutive arcs of some cycle~$ \gamma $.
	So $ (e_3 \cup e_1) $ and $ (e'_1 \cup e'_3) $ are two consecutive
	arcs of some cycle $ \gamma' $.  
    
\end{enumerate}
Notice that in each case $ \gamma' $ may be the same cycle than $ \gamma $.

\begin{defn}
We will say that two thick-Y's are \emph{glued with a twist} in
the case (2) above, and \emph{glued without a twist} in the case (1) above
(see fig.~\ref{figure10}).
\end{defn}

\begin{figure}[ht]
\includegraphics{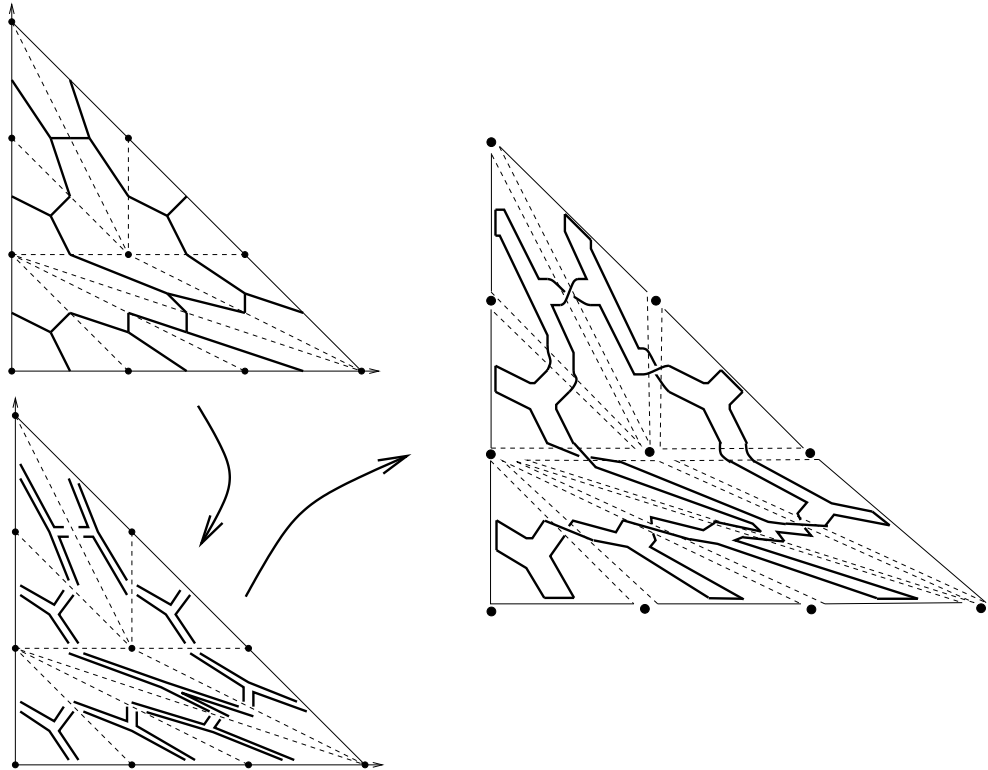}  
\caption{\small To construct $ F(K) $:  Thicken the Y's in each triangle and
glue any two adjacent thick-Y's with or without a twist (depending on $ K
$).}
\label{figure11}  
\end{figure}

Since $ G(\Pi) $ is connected, we obtain, by gluing in this way
the thick-Y's of all the triangles of $ \cT $, a surface with boundary (see
fig.~\ref{figure11}).

\begin{defn}
The surface with boundary constructed above from a T-curve $ K $ will
be called a \emph{T-filling}, and will be denoted $ F(K) $.
\end{defn}

Notice from the construction of the T-filling $ F(K) $, that the connected
components of $ \partial F(K) $ are in one-to-one correspondence with the
connected components of $ K $.

\subsection{The Relation with Algebraic Geometry} 
\label{F(K)-alg}

The T-filling of a T-curve is in fact a topological analog of the quotient of
the complexification of an algebraic curve by the complex conjugation.  Since
the complexification of an algebraic curve is better known for curves on the
projective plane than for curves on arbitrary projective toric surfaces, we
restrict here the comparison to that case, bearing in mind that it holds in
the more general case.

A real plane projective curve $ C = C(f) $ of degree $ d $ is the subset $
\{(x_0 \colon x_1 \colon x_2) \in \RP^2, f(x_0, x_1, x_2) = 0\} $ for some
homogeneous polynomial $ f $ of degree $ d $ with real coefficients.  The
complexification $ \C C $ of $ C $ is the surface $ \{(z_0 \colon z_1 \colon
z_2) \in \CP^2, f(z_0, z_1, z_2) = 0\} $ in $ \CP^2 $.  The curve is
non-(complex)-singular if the derivative is a nonzero function in all points
of the curve (of its complexification).  Notice that the complex conjugation
leaves $ \C C $ invariant, and that $ C $ is its set of fixed points.
Therefore the quotient of $ \C C $ by the complex conjugation is a surface
with $ C $ as boundary.

Let $ K = K(\gT_d, \cT, \delta) $ be a T-curve on $ \RP^2 $ = $ S(\gT_d) $
(see remark~\ref{amb-tor}).  We have seen in remark~\ref{Tcurv/tor} that ,
when the triangulation $ \cT $ of $ T_d $ is convex, $ K $ correspond to a
real plane projective non singular curve $ C $ of degree $ d $.  

\begin{figure}[h]
\includegraphics{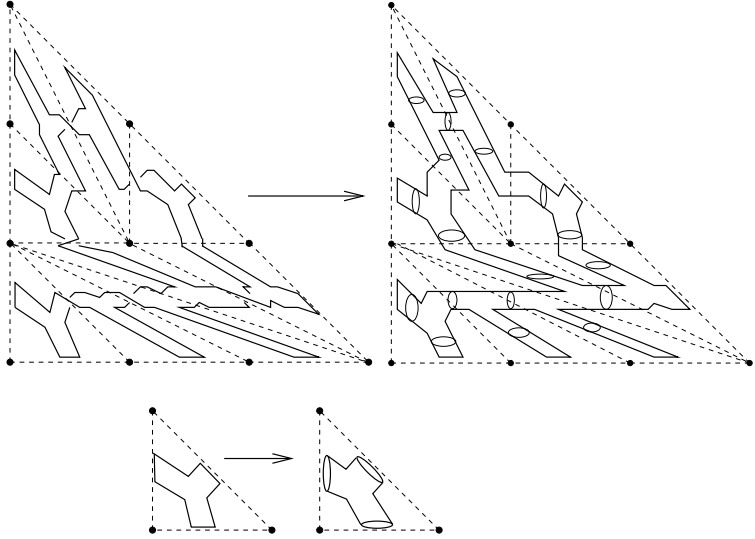}
\caption{\small The double of this T-curve on $ \RP^2 = S(\gT_3) $ is a
surface of genus $ \frac{(3-1)(3-2)}{2} = 1 $.}
\label{figure12}
\end{figure}

It is well known that $ \C C $ is an orientable surface of genus $
\frac{(d-1)(d-2)}{2} $.  Similarly the double of $ F(K) $ along its boundary
is a surface of genus $ \frac{(d-1)(d-2)}{2} $.  Indeed to double it we
just have to let the (flat) thick-Y's become ``hollow-Y's'', glue them
without worrying about the twists, and close the open ends on $ \bd\gT $ by
disks (see fig.~\ref{figure12}).  The number of handles of this surface is
equal to the number of interior integer points of $ \gT $, and this is
precisely $ \frac{(d-1)(d-2)}{2} $ (see fig.~\ref{figure12}).  Moreover the
number of connected components of the boundary is the number of connected
components of $ K $.  Therefore $ F(K) $ is homeomorphic to the quotient of $
\C C $ by the complex conjugation.

The analogy goes further.  The surface quotient of $ \C C $ by the complex
conjugation allows to define the orientability and, if it is orientable, an
orientation of $ C $.  Similarly a notion of orientability and orientation of
$ K $ emerges naturally from $ F(K) $ (see section~\ref{sec:orient-T}).

\begin{rem}
\label{Tfill-high-dim}
If the T-filling of a T-surface is to be homeomorphic to the quotient of the
complexification of a real algebraic surface, it is of dimension four.  It is
not clear how to construct this four dimensional variety in a similar way
than the T-filling of a T-curve.
\end{rem}

\subsection{Applications}
\label{Tfill-applic}

\subsubsection{First Application: The Harnack Theorem for T-Curves}
\label{sec:Har-T}

\begin{thm}[``Harnack'' for T-curves] 
\label{th:Har-T}
\one\ The number of connected components of a T-curve is no more than
the number of integer points in the interior of its Newton polygon, plus
one. \\
\two\ There are T-curves with arbitrary Newton polygons which
achieve the upper bound above.
\end{thm}

\begin{pf}
Let us prove just \one\ for the moment.  \two\ will be proved in
prop.~\ref{Har-cocp} where we will construct T-curves with arbitrary
Newton polygon which achieve the upper bound.

Let $ K(\Pi, \cT, \delta) $ be a T-curve, and let $ i $ be the number
of integer points in the interior of $ \Pi $.  Let us glue disks along $
\bd F(K) $ to obtain from $ F(K) $ a closed surface $ \Sigma = \Sigma(K) $.
It is clear from the construction of $ F(K) $ that the incidence graph $
G(\Pi) $ of the triangulation $ \cT $ is a retraction of $ F(K) $, so let us
consider $ \Sigma $ as the cell complex obtained by gluing disks on $ G(\Pi)
$ along the cycles $ \mu(\gamma) $.  

Let $ D $ be the number of disks glued to $ G(\Pi) $, and let $ E' $ and $ V'
$ be respectively the number of edges and vertices of $ G(\Pi) $.  Let $ T $,
$ E $ and $ V $ be respectively the numbers of triangles, edges and vertices
of the triangulation $ \cT $ of $ \Pi $, and let $ L $ be the length of $ \bd
\Pi $.  Each triangle of $ \cT $ contains three edges of $ G(\Pi) $ so $ E' =
3T $.  The vertices of $ G(\Pi) $ lie in each triangle and each edge of $ \cT
$, so $ V' = T+E $.

On one hand the Euler characteristic of $ \Sigma $ is 

\begin{eqnarray*} 
\chi(\Sigma) 
& = & D - E' + V'  \\
& = & D - 3T + (T+E)  \\
& = & D + T - E + L \qquad\text{$ \cT $ is a triangulation, therefore 
$ 3T = 2E - L $} \\ 
& = & D + 1 - V + L \qquad\text{Euler relation on $\cT$ implies 
$ T - E + V = 1 $} \\
\end{eqnarray*}

On the other hand $ \chi(\Sigma) $ is equal to $ 2 - 2g $ where $ g $ is the
genus of $ \Sigma $.  Therefore we get

\begin{eqnarray*}
D & = &  (V - L) + 1 - 2g  \\
& \leq & i + 1 \qquad\text{since $ V - L = i $ and $ g \geq 0 $}
\end{eqnarray*}
This proves \one\ since $ D $ is also the number of connected components of $
K $.  
\qed
\end{pf}

\begin{defn}
A T-curve $ K(\Pi, \cT, \delta) $ with the number of connected components
equal to the number of integer points in the interior of $ \Pi $ plus one
is called a \emph{maximal} T-curve.  
\end{defn}

It is well known that a real projective curve is maximal if and only if the
quotient of its complexification by the complex conjugation is homeomorphic
to a sphere with a positive number of points removed.  In particular it has
no handle.  We have here the analog statement in the case of T-curves:

\begin{cor}
\label{max-A-disk}
A T-curve is maximal if and only if its T-filling is homeomorphic to a sphere
with a positive number of points removed.  
\end{cor}

\begin{pf}
This statement follows from the proof of~\ref{th:Har-T}.  Indeed the
curve is maximal if the surface $ \Sigma(K) $ has genus $ 0 $.  Since a
T-curve has at least one connected component, $ F(K) $ is obtained from $
\Sigma(K) $ by removing a positive number of disks.
\qed
\end{pf}

\subsubsection{Harnack T-Curves}
\label{Har-const}

Let $ K $ be a T-curve with Newton polygon $ \Pi $.  Recall that the
quadrants of the ambient surface $ S(\Pi) $ are represented in $ \R^2 $ by
the copies $ (\sig_{a,b} \cdot \Pi) = \{(x,y), ((-1)^a x, (-1)^b y) \in \Pi\}
$, where $ a,b \in \{0,1\} $.

\begin{defn}
\label{df:Har-dist}
Let $ c,a,b \in \{0,1\} $, and let $ \delta $ be the distribution of signs
defined on the points $ (x,y) \in (\sig_{a,b} \cdot \Pi) $ by: 
$$ 
\delta(x,y) = \left\{ \begin{array}{ll}
		    (-1)^c & \text{if $ (x,y) $ is of even parity} \\
		    (-1)^{c+1} & \text{if $ (x,y) $ is of odd  parity} 
		    \end{array} \right. 
$$
The distribution of sign of a T-curve deduced from $ \delta $ by the
formula~\ref{ext-sd2} page~\pageref{ext-sd2} will be called a \emph{Harnack
sign distribution}, and $ (c,a,b) $ will be called the \emph{type} of $
\delta $.  A T-curve with a Harnack sign distribution will be called a
\emph{Harnack T-curve}.
\end{defn}

In order to write a unique formula to describe a Harnack distribution of
sign, let $ [P] $ be equal to $ 1 $ if the proposition $ P $ is true, and
equal to $ 0 $ if the proposition $ P $ is false.  Now let $ g,h = 0 $ or $ 1
$, let $ (x,y) $ be an integer point of $ (\sig_{g,h} \cdot \Pi) $, let $
(e,f) $ be the parity of $ (x,y) $, and let $ (c,a,b) $ be the type of a
Harnack distribution of signs $ \delta $ on $ \Pi $.  Then the formula of the
above definition (\ref{df:Har-dist}), is rewritten:
\begin{equation}
\label{eq:Har-sign}
\delta(x,y) = (-1)^{c + [(e,f) \not = (0,0)] + \langle (e,f), (g+a, h+b)
\rangle}
\end{equation}

\begin{defn}
Let $ \Pi $ be a polygon, and let $ a $ be an arc which splits $
\Pi $ into two connected components.  We will say that $ a $ \emph{surrounds
an integer point $ P $ in $ \Pi $} if $
P $ is the only integer point in the closure of one of the connected
component of $ \Pi \setminus a $.
\end{defn}

\begin{lem}
\label{pts-surrounded}
Let $ K(\Pi, \cT, \delta) $ be a Harnack curve, and let $ (c,a,b) $ be the
type of $ \delta $.  Then for every $ g,h = 0 $ or $ 1 $, the integer points
of parity $ (b+h, a+g) $ will be surrounded
in the quadrant $ (\sig_{g,h} \cdot \Pi)
$ by an oval of $ K $ if they belong to the interior of the quadrant, or by
an arc of $ K $ if they are on the boundary of the quadrant.
\end{lem}

\begin{pf}
According to formula~\ref{eq:Har-sign} above, the sign of an integer
point of parity $ (e,f) $ in quadrant $ (g,h) $ depend only on the value of
the expression
\begin{equation}
\label{eq:log-sign}
[(e,f) \not = (0,0)] + \langle (e,f), (g+a,h+b) \rangle
\end{equation}
The following array displays the value of this expression for each parity $
(e,f) $ (one parity per line) in each quadrant $ (g+a,h+b) $ (one quadrant
per column):
$$
\begin{array}{c@{\quad}||@{\quad}c@{\quad}|@{\quad}c@{\quad}|@{\quad}%
c@{\quad}|@{\quad}c@{\quad}|}
      & (0,0) & (1,0) & (1,1) & (0,1) \\
\hline \hline
(0,0) & 0+0=0 & 0+0=0 & 0+0=0 & 0+0=0 \\
\hline
(1,0) & 1+0=1 & 1+1=0 & 1+1=0 & 1+0=1 \\
\hline
(1,1) & 1+0=1 & 1+1=0 & 1+0=1 & 1+1=0 \\
\hline
(0,1) & 1+0=1 & 1+0=1 & 1+1=0 & 1+1=0 \\
\hline
\end{array}
$$
By looking at each row of the array we see that in a given quadrant, points
of parity $ (e,f) = (b+h, a+g) $ have opposite signs than their neighbors.
This implies that they are separated from their neighbors by $ K $ (by an arc
if the point is on the boundary of the quadrant, and by an oval otherwise). 
\qed
\end{pf}

\begin{figure}[ht]
\includegraphics{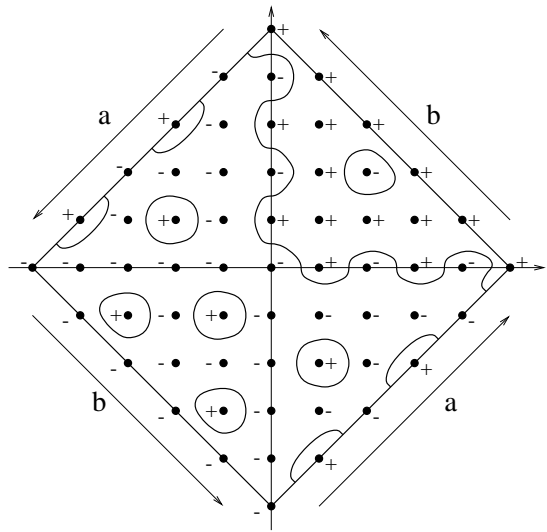}  
\caption{\small A Harnack T-curve of degree $ 5 $.}
\label{figure13}  
\end{figure}

\begin{figure}[ht]
\includegraphics{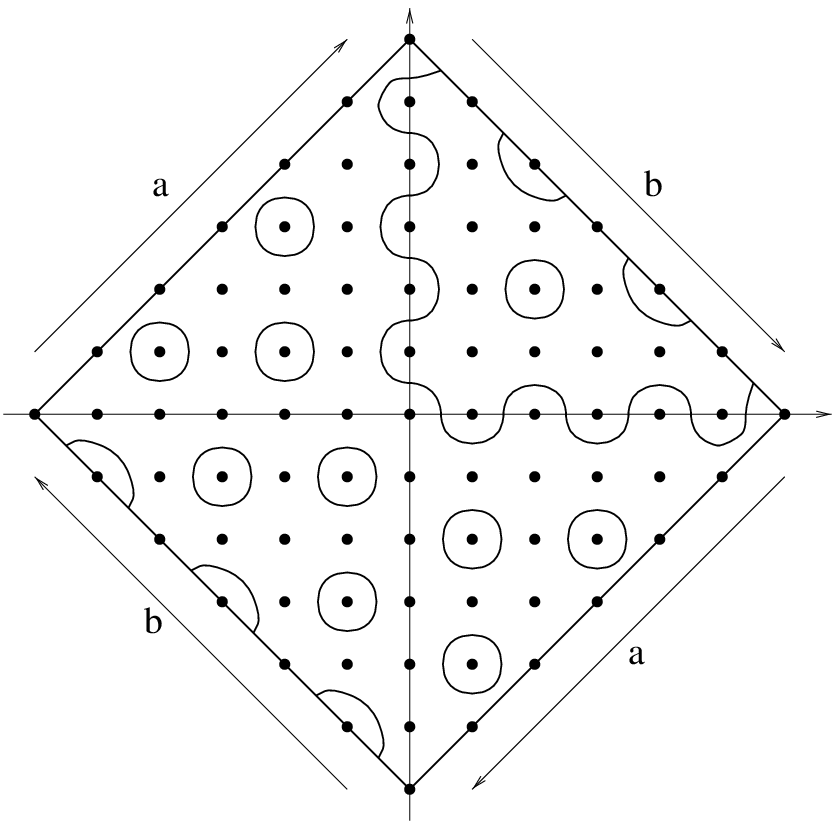}  
\caption{\small A Harnack T-curve of degree $ 6 $.}
\label{figure14}  
\end{figure}

Harnack T-curves are called so because on $ \RP^2 $ they are congruent to the
well known curves constructed by Harnack~\cite{Harnack1} which have 
\begin{itemize}

\item  One one-sided component and $ \frac{(d-1)(d-2)}{2} $ outermost empty
	ovals if they have odd degree $ d $ (see example fig.~\ref{figure13}).

\item  One outermost oval containing $ \frac{(k-1)(k-2)}{2} $ empty ovals
	and $ \frac{3k(k-1)}{2} $ outermost empty ovals if they have even
	degree $ d = 2k $ (see example fig.~\ref{figure14}).

\end{itemize}

The additive group $ \Ztwo^3 = \{0,1\} \times \{0,1\} \times \{0,1\} $ acts
on the distributions of signs as follows:  Let $ \theta = (c,a,b) \in \Ztwo^3
$, and let $ \delta $ be a distribution of signs on some integer points.
Then 
\begin{equation} 
\label{theta-sd}
(\theta \cdot \delta)(x,y) = (-1)^{c + \langle (a,b), (x,y) \rangle}
\delta(x,y) 
\end{equation}
Therefore $ \Ztwo^3 $ acts also on the set of T-curves with a given
Newton polygon:
$$
\theta \cdot K(\Pi, \cT, \delta) = K(\Pi, \cT, \theta \cdot \delta) 
$$
Observe that $ \theta \cdot \delta(x,y) = (-1)^c \delta(\sig_{a,b} \cdot
(x,y)) $, and since $ K(\Pi, \cT, \delta) = K(\Pi, \cT, -\delta) $, the group
$ \Ztwo^3 $ must be considered as $ \Ztwo \times (\Ztwo^2) $, where the first
factor $ \Ztwo $ acts trivially and the second factor $ \Ztwo^2 = \{(a,b),\:
a,b = 0 \,\text{or}\, 1\} $ acts as the group of symmetries $ \{\sig_{a,b}:
(x,y) \mapsto ((-1)^a x, (-1)^b y)\} $ on the four quadrants of
$ S(\Pi) $.

\begin{lem}
\label{ref-Har-dist}
A Harnack T-curve $ K'(\Pi, \cT, \delta') $, with $ \delta' $ of type $
\theta' = (c,a,b) $, is the image by the symmetry $ \sig_{a,b} $ of the
Harnack T-curve $ K(\Pi, \cT, \delta) $ with $ \delta $ of type $ (1,0,0) $.
\end{lem}

\begin{figure}[ht]
\includegraphics{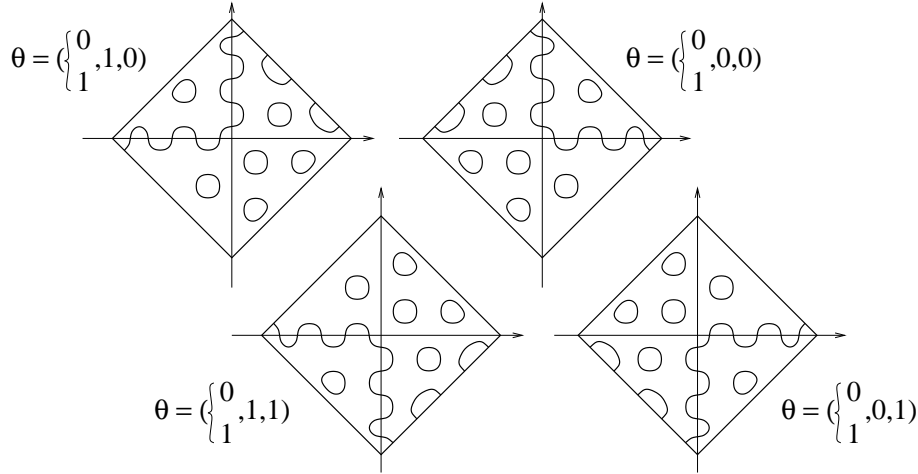}  
\caption{\small Eight types of Harnack distributions for four Harnack curves
symmetric to one another.}
\label{figure15}  
\end{figure}

\begin{pf}
Indeed from the definition~\ref{df:Har-dist} of the Harnack
distributions and from formula~\ref{theta-sd} above we get that, for any $
\theta \in \Ztwo^3 $, the distribution $ \theta \cdot \delta' $ is a Harnack
distribution of type the sum $ \theta + \theta' $.  So for $ \theta =
(c+1,a,b) $ we get that $ \theta \cdot K' = K $.  This is equivalent to $ K'
= \theta \cdot K $.  Since $ \theta \cdot K = K(\Pi, \cT, \theta \cdot
\delta) $ we get from the observation above that $ K' = \sig_{a,b}(K) $.
see fig.~\ref{figure15}.
\qed
\end{pf}

\begin{prop}
\label{Har-cocp}
Let $ K = K(\Pi, \cT, \delta) $ be a Harnack T-curve with Harnack
distribution of type $ (c,a,b) $.  Let $ g(s,t) $ be the number of points of
parity $ (s,t) $ in the interior of $ \Pi $.  Then 
\begin{enumerate}

\item  \label{Har-max}
	$ K $ is a maximal T-curve. 

\item  The connected components of $ K $ are distributed on its ambient
	surface as follows:

\begin{enumerate}

\item  \label{cocp-ev}
	$ g(0,0) $ empty ovals of sign\footnote[2]{recall the definition of the
	sign of an oval is def.~\ref{df:sign-ov}} $ (-1)^c $ in the quadrant
	$ (\sig_{a,b} \cdot \Pi) $.

\item  \label{cocp-odd} 
	$ g(s,t) $ empty ovals of sign $ (-1)^{c+1} $ in the
	quadrant $ \sig_{t,s} \cdot \sig_{a,b} \cdot \Pi = \sig_{t+a,s+b} 
	\cdot \Pi $, for each odd parity $ (s,t) $.

\item  \label{last-ntr}
	either a nontrivial connected component if $ \Pi $ has one broken
	edge of odd length.

\item  \label{last-ov}
	or an oval surrounding exactly all the ovals of \textit{(a)} if
	all the broken edges of $ \Pi $ are of even length.

\end{enumerate}
\end{enumerate}
\end{prop}

\begin{pf}
Thanks to lemma~\ref{ref-Har-dist} we assume without loss of
generality that $ (c,a,b) = (1,0,0) $.  

From lemma~\ref{pts-surrounded} we get that every point of even parity which
lies in the interior of $ \Pi $ is surrounded
 by an oval of $ K $.  From the
definition~\ref{df:Har-dist} of Harnack distribution of sign, we know that
this oval is of sign $ -1 $.  This proves assertion~\ref{cocp-ev}.

From lemma~\ref{pts-surrounded} we get as well that for every point $ P $ of
odd parity $ (s,t) $ which lies in the interior of $ \Pi $, its copy $
(\sig_{t,s} \cdot P) $ is surrounded by an oval of $ K $.  From the
definition~\ref{df:Har-dist} of Harnack distribution, and from
formula~\ref{eq:Har-sign} page~\pageref{eq:Har-sign} we get that this oval is
of sign $ +1 $.  This proves assertion~\ref{cocp-odd}.

Let $ g = \sum_{s,t} g(s,t) $.  We just proved that $ K $ has $ g $ empty
ovals.  Moreover from lemma~\ref{pts-surrounded} we get that for every
integer point $ \bd \Pi $ one of its symmetric copy is surrounded
by an arc of $ K $, so there is at least one more connected component of $ K
$.  From the part \one\ of theorem~\ref{th:Har-T} we deduce that $ K $ has
exactly $ g+1 $ connected component.  This proves assertion~\ref{Har-max}.

Let us call $ O $ the connected component of $ K $ which intersects
(transversally) all the lifts of the broken edges of $ \Pi $.

Assume first that $ \Pi $ has a broken edge $ \ell $ of odd length.  From
lemma~\ref{2-e-neg} we know that among the two lifts of an edge $ e $ of the
triangulation $ \cT $ lying on $ \ell $ one is of negative sign and the other
one is of positive sign.  Therefore the lift $ \mu^{-1}(\ell) $ in $ S(\Pi) $ is
intersected transversally an odd number of times by $ O $, which implies,
since $ \mu^{-1}(\ell) $ is an embedded circle, that $ O $ is not an oval of $ K
$ (see fig.~\ref{figure16}).  This proves assertion~\ref{last-ntr}.

\begin{figure}[ht]
\includegraphics{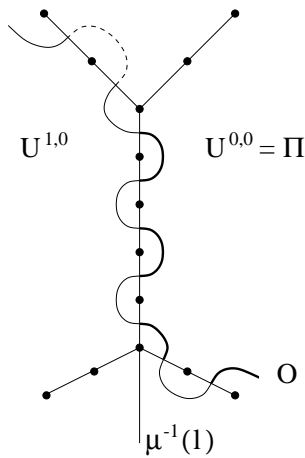}  
\caption{\small If a broken edge of $ \Pi $ is of odd length, then the
component $ O $ of $ K $ cuts the embedded circle $ \mu^{-1}(\ell) $ an odd
number of times.}
\label{figure16}  
\end{figure}

Assume now that all the broken edges of $ \Pi $ are of even length.  This
implies that all the endpoints $ \mu^{-1}(u_i) $ of the broken edges are of
same parity.  Thanks to prop.~\ref{translat-Pi}, we assume without loss of
generality that they are of even parity.  Therefore one integer point over
two on $ \bd \Pi $ is of even parity.  Since $ \delta $ is assumed to be of
type $ (1,0,0) $, all the points of even parity lying on $ \bd \Pi $ are
surrounded by arcs of $ \Pi $ in the quadrant $ \Pi $, and all the other
points lying on $ \bd \Pi $ are surrounded in another quadrant.

By retracting each arc, which surrounds an integer point lying on $ \bd \Pi
$, toward this point, we shrink $ O $ onto the boundary of the quadrant $ \Pi
$.  Since $ \Pi $ is homeomorphic to a disk, this boundary, and hence $ O $,
is an oval.  It is clear that during the shrinking no crossing with another
connected component of $ K $ happened, so $ O $, like the boundary of $ \Pi
$, surrounds exactly the empty ovals lying within $ \Pi $ (see
fig.\ref{figure17}). This proves assertion~\ref{last-ov}.
\qed
\end{pf}

Notice that the assertion~\ref{Har-max} of proposition~\ref{Har-cocp} proves
\two\ of theorem~\ref{th:Har-T}, and therefore finishes the proof
of~\ref{th:Har-T}.
\qed

\begin{figure}[ht]
\includegraphics{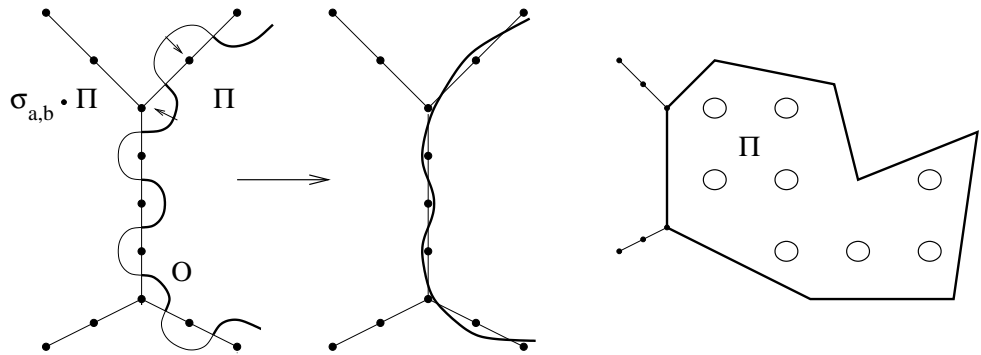}  
\caption{\small When all the edges of $ \Pi $ are of even length, the
component $ O $ of $ K $ can be shrank onto the boundary of $ \Pi $.} 
\label{figure17}  
\end{figure}

\begin{cor}
\label{Har-cong}
Any two Harnack T-curves with same Newton polygon are congruent by an
homeomorphism of $ S(\Pi) $ which is the identity on the boundary of the
quadrants of $ S(\Pi) $.
\end{cor}

\begin{pf}
Thanks to lemma~\ref{ref-Har-dist} we restate without loss of generality this
corollary for any two Harnack T-curves with sign distribution of type $
(1,0,0) $.  Now since the sign distributions on the integer points are the
same, the sign of the segments on the boundary of the quadrants of $ S(\Pi) $
are also the same.  This fixes the T-curves on this boundary.
Prop.~\ref{Har-cocp} shows that the connected components of any two such
Harnack T-curves, restricted to a given quadrant, are congruent. 
\qed
\end{pf}

\subsubsection{Second Application: Orientation of T-Curves}
\label{sec:orient-T}

Let $ C $ be a real plane projective nonsingular curve of degree $ d $, let $
\C C $ be its complexification, and let $ \tau $ be the complex conjugation.
Recall that $ \C C $ is an orientable surface of genus $ \frac{(d-1)(d-2)}{2}
$ and that $ C $ is the set of fixed points of $ \tau $.   

\begin{defn}
\label{alg-type-I-II}
$ C $ is called a \emph{dividing curve} (or a curve of
\emph{type I}) if it divides its complexification into two connected
components, and is called a \emph{non-dividing curve} (or a curve of
\emph{type II}) if it does not divide its complexification.
\end{defn}

\begin{lem}
$ C $ is of type~I if and only if the quotient of its complexification by the
complex conjugation is orientable.  
\end{lem}

\begin{pf}
If $ C $ is of type I, then the quotient $ \C C / \tau $ is homeomorphic to a
connected component of $ \C C \setminus C $, and therefore is orientable.

Assume now $ C $ is of type II\@.  Let $ P \in \C C $ a point not in $ C $,
and let $ \gamma $ be a path connecting $ P $ to its conjugate $ \tau(P) $
which doesn't intersect $ C $.  Since $ \tau $ reverses the orientation of $
\C C $, the image of the path $ \gamma $ in $ \C C / \tau $ is a closed path
which reverses the orientation.
\qed
\end{pf}

By analogy we can now define a notion of orientability on T-curves via the
T-filling in the following way. 

\begin{defn}
\label{Tcurv-type-I-II}
A T-curve is of \emph{type I} if its T-filling is orientable, and of
\emph{type II} if its T-filling is non-orientable.
\end{defn}

Let $ C $ be a real plane projective nonsingular curve of type I\@.  The
curve $ C $ divides its complexification $ \C C $ into two orientable halves.
An orientation on $ \C C $ induces an orientation on each of the halves.
These orientations induce in turn two opposite orientations on $ C $.

\begin{figure}[h]
\includegraphics{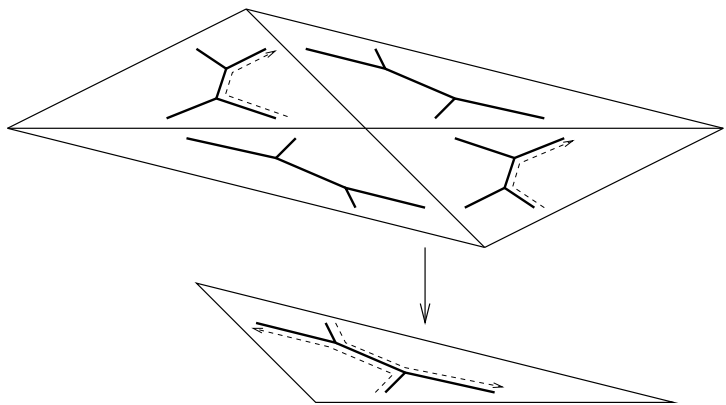} 
\caption{\small The oriented cycles on the incidence graph on $ \Pi $ lift to
oriented cycles on the incidence graph on $ S(\Pi) $.}
\label{figure18} 
\end{figure}

By analogy we can now define a notion of orientation on any T-curve $ K =
K(\Pi, \cT, \delta) $ of type I via its T-filling $ F(K) $ in the following
way.  Since $ F(K) $ is orientable, an orientation on it induces an
orientation on its boundary.  By letting $ F(K) $ retract to the incidence
graph $ G(\Pi) $, we get an orientation of the cycles $ \gamma =
\mu(\tilde{\gamma}) $ which are the images of the connected components of $ K
$.  An segment $ e $ of $ \gamma $ lifts to a segment $ \sig_{a,b} \cdot e $
of $ \tilde{\gamma} $.  If $ e $ is oriented from endpoint $ (x,y) $ to
endpoint $ (x',y') $, then segment $ \sig_{a,b} \cdot e $ will be oriented
from $ \sig_{a,b} \cdot (x,y) $ to $ \sig_{a,b} \cdot (x',y') $.  We get this
way an orientation of all the connected components $ \tilde{\gamma} $ of $ K
$ (see fig.~\ref{figure18}).  The two opposite orientations on $ F(K) $
induce then two opposite orientations on $ K $.

\begin{defn}
\label{I-orient}
The orientation described above on a T-curve of type I will be called
a \emph{(type I)-orientation} of the T-curve.
\end{defn}

Now we have defined type I and type II, and type-I-orientation for a T-curve,
it would be interesting to prove for T-curves of given degree on $ \RP^2 $
theorems known for real plane projective nonsingular curves which take into
account type and orientation (in particular Arnol'd
congruence~\cite{Arnold/ovals-invol-4mfds-quad-forms}, and Rokhlin
formula~\cite{Rokhlin/cplex-top-char}).  It would be a first step in order to
find extensions of these theorems to T-curves on arbitrary ambient surfaces
and to real algebraic curves on arbitrary toric varieties.

\bibliography{articles}
\bibliographystyle{acm}

\end{document}